\documentclass[a4paper,11pt]{amsart}

\def\frak{\mathfrak}
\def\Bbb{\mathbb}
\def\Cal{\mathcal}

\newtheorem*{prop*}{Proposition}

\newtheorem*{thm*}{Theorem}

\newtheorem*{lem*}{Lemma}

\newtheorem*{kor*}{Corollary}

\newcommand{\id}{\operatorname{id}}

\newcommand{\vol}{\operatorname{vol}}

\newcommand{\x}{\times}
\renewcommand{\o}{\circ}

\let\ccdot\cdot
\def\cdot{\hbox to 2.5pt{\hss$\ccdot$\hss}}

\newcommand{\al}{\alpha}
\newcommand{\be}{\beta}
\newcommand{\ga}{\gamma}

\newcommand{\ep}{\varepsilon}
\newcommand{\ka}{\kappa}

\newcommand{\om}{\omega}
\renewcommand{\phi}{\varphi}
\newcommand{\ph}{\varphi}

\newcommand{\si}{\sigma}

\newcommand{\La}{\Lambda}
\newcommand{\Ga}{\Gamma}
\newcommand{\Ph}{\Phi}
\newcommand{\Om}{\Omega}

\newcommand{\ce}{{\Cal E}}

\usepackage{amssymb}
\usepackage{amscd}

\newcommand{\nd}{\nabla}
\newcommand{\Ps}{\Psi}
\newcommand{\Rho}{{\mbox{\sf P}}}

\newcommand{\End}{\operatorname{End}}

\newcommand{\Ric}{\operatorname{Ric}}

\def\sideremark#1{\ifvmode\leavevmode\fi\vadjust{\vbox to0pt{\vss
 \hbox to 0pt{\hskip\hsize\hskip1em
 \vbox{\hsize3cm\tiny\raggedright\pretolerance10000
 \noindent #1\hfill}\hss}\vbox to8pt{\vfil}\vss}}}%
\newcommand{\cq}{{\Cal Q}}                     
\newcommand{\aM}{\tilde{M}}                    

\newcommand{\X}{\mbox{\boldmath{$ X$}}}        
\newcommand{\sX}{\mbox{\scriptsize\boldmath{$X$}}}        
\newcommand{\h}{\mbox{\boldmath{$ h$}}}        

\newcommand{\aR}{\mbox{\boldmath{$ R$}}}             
\newcommand{\aRic}{\mbox{\boldmath{$Ric$}}}
\newcommand{\annd}{\mbox{\boldmath$ \nabla$}}  
\newcommand{\cce}{\tilde{\ce}}                
\newcommand{\adD}{\mbox{{$\tilde{\mathbb{D}}$}}}       
\newcommand{\dD}{\mbox{{$\mathbb{D}$}}}       
\newcommand{\D}{\mbox{\boldmath{$ D$}}}       
\newcommand{\aDe}{\mbox{\boldmath$ \Delta$}}  
\newcommand{\V}{{\mbox{\boldmath$ V$}}}

\newcommand{\tf}{\tilde{f}}

\begin{document}
\title[tractors and conformal ambient metric]{Standard Tractors and
  the Conformal Ambient Metric Construction}  
\author{Andreas \v Cap and A.\ Rod Gover}
\date{July 1, 2002}
\address{A.C.: Institut f\"ur Mathematik, Universit\"at Wien, Strudlhofgasse 4,
A--1090 Wien, Austria and International Erwin Schr\"odinger Institute for
Mathematical Physics, Boltzmanngasse 9, A--1090 Wien, Austria\newline\indent
A.R.G.: Department of Mathematics, The University of Auckland, Auckland, 
New Zealand}
\email{Andreas.Cap@esi.ac.at, gover@math.auckland.ac.nz}
\begin{abstract}
In this paper we relate the Fefferman--Graham ambient
metric construction for conformal manifolds to the approach to
conformal geometry via the canonical Cartan connection. We show that
from any ambient metric that satisfies a weakening of the usual
normalisation condition, one can construct the conformal standard
tractor bundle and the normal standard tractor connection, which are
equivalent to the Cartan bundle and the Cartan connection. This result
is applied to obtain a procedure to get tractor formulae for all
conformal invariants that can be obtained from the ambient metric
construction. We also get information on ambient metrics which
are Ricci flat to higher order than guaranteed by the results of
Fefferman--Graham.
\end{abstract}

\maketitle

\section{Introduction}\label{1}
It is an old result of E.~Cartan that conformal manifolds of dimension
$\geq 3$ admit a canonical normal Cartan connection. While this solves
the equivalence problem for conformal structures, the problem of
a complete conformal invariant theory and the related problem of
constructing conformally invariant differential operators remain very
difficult. Based on similar ideas for CR--structures, Ch.~Fefferman
and C.R.~Graham initiated a new approach to these problems in 1985,
see \cite{F-G}. Viewing a pseudo--conformal structure of signature
$(p,q)$ on a manifold $M$ as a ray bundle $S^2T^*M\supset\cq\to M$ the
idea of this approach is to associate to the given conformal structure
a pseudo--Riemannian metric of signature $(p+1,q+1)$ on $\cq\x
(-1,1)$, the so--called ambient metric. This metric is required to be
homogeneous and compatible with the conformal structure in a rather
obvious sense, while the more subtle condition is that it should be
Ricci flat up to some order along $\cq$. In odd dimensions, ambient
metrics which are Ricci flat to infinite order along $\cq$ exist and
are essentially unique, while in even dimensions there is an
obstruction at finite order, but up to that order the ambient metric
is again essentially unique. 

This immediately leads to a construction for conformal invariants,
since any Riemannian invariant of an ambient metric satisfying these
Ricci conditions, which is of low enough order in the
even--dimensional case, gives rise to a conformal invariant of the
underlying conformal structure. Moreover, the ambient metric
construction has been applied in \cite{GJMS} to construct conformally
invariant powers of the Laplacian. This constructions gives
arbitrarily high powers in odd dimensions and powers up to some
critical order in even dimensions.  This is complemented by
C.R.~Graham's result in \cite{Graham} that there is no conformally
invariant third power of the Laplacian in dimension four, which
strongly suggests that the obstruction to the ambient metric
construction in even dimensions is of fundamental nature. It should
also be remarked here that the ambient metric construction has
recently received renewed interest because of its relation to the
so--called Poincar\'e metric and via that to scattering theory and the
AdS/CFT--correspondence in physics, see \cite{W,GrW,HS2,GrZ,GrF} and
references therein.

Over the past few years the Cartan approach to conformal geometry and
a more general class of geometric structures called parabolic
geometries has been significantly developed. One development, whose
origins can be traced back to the work of T.~Thomas in the 1920's and
1930's \cite{T,Thomasbook}, is the concept of tractor bundles, which
give an equivalent description of the Cartan bundle and the Cartan
connection in terms of linear connections on certain vector bundles.
These then lead to an efficient calculus, which has been successfully
applied to the study of conformal invariants and conformally invariant
differential operators, see e.g.~\cite{BEG, Mike-Srni, conf-invar}.

The purpose of this paper is to relate precisely the ambient metric
construction to the conformal standard tractor bundle and its
canonical linear connection. We first construct a standard tractor
bundle and a tractor connection on that bundle from a very general
class of ambient metrics. Then we prove that normality of this tractor
connection is equivalent to vanishing of the tangential components of
the Ricci curvature of the ambient metric along $\cq$. Hence, we
obtain the normal standard tractor bundle and tractor connection from
any metric produced by the ambient metric construction. This is done
in section \ref{2}.

In section \ref{3}, we express some basic elements of tractor calculus
in terms of ambient data. This is then used to show that for any
ambient metric which satisfies the Ricci conditions of \cite{F-G},
there is an algorithm to compute all covariant derivatives of the
curvature (up to the critical order in even dimensions) from the
tractor curvature. Since any local scalar conformal invariant obtained
from the ambient metric construction is a complete contraction of a
tensor product of such covariant derivatives, we obtain an algorithm
to compute a tractor formula for any of these invariants. This is an
important achievement, since in contrast to the situation of the
ambient metric construction, converting tractor formulae into formulae
in terms of metrics from the conformal class is a purely mechanical
procedure and in particular does not involve solving any equations. In
fact it is straightforward to write software for these expansions, see
\cite{Gover-Peterson}.

Our results also cover the case of metrics which are Ricci flat to
higher order than the ones that can be obtained from the ambient
metric construction. While such metrics do not exist on general
conformal manifolds, studying the cases in which they do exist is of
considerable interest in conformal geometry. What we can prove in this
case is that all higher covariant derivatives of the curvature can be
obtained from the tractor curvature and one ``critical'' covariant
derivative. 

It should be pointed out that formally, our results are completely
independent from the results of \cite{F-G} on the ambient metric
construction. In particular, Theorem \ref{3.4} goes a long way towards
an independent proof of the uniqueness of the ambient metric. In fact, it
contains enough information on the uniqueness to show that the ambient
metric construction leads to conformal invariants. We believe that our
results can be extended to a complete proof of uniqueness of the
ambient metric, which is of entirely different nature than the one in
\cite{F-G}. The reason why we do not go further in that direction here
is that we believe that the ideas we develop also can be used for
existence proofs for ambient metrics and an analysis for the
obstruction to the existence of an ambient metric which is Ricci flat
to infinite order in even dimensions, and we will  take up this
whole circle of problems elsewhere. 

The authors would like to thank C.\ Robin Graham for several helpful
conversations.

\section{The ambient construction of conformal standard
tractors}\label{2} 
\subsection{Conformal structures}\label{2.1}
Let $M$ be a smooth manifold of dimension $n\geq 3$. A {\em conformal
structure\/} on $M$ of signature $(p,q)$ (with $p+q=n$) is an
equivalence class of smooth pseudo--Riemannian metrics of signature
$(p,q)$ on $M$, with two metrics being equivalent if and only if one
is obtained from the other by multiplication with a positive smooth
function.

For a point $x\in M$, and two metrics $g$ and $\hat g$ from the
conformal class, there is an element $s\in\Bbb R_+$ such that $\hat
g_x=sg_x$. Thus, we may equivalently view the conformal class as being
given by a smooth ray subbundle $\cq\subset S^2T^*M$, whose fibre at $x$ is
formed by the values of $g_x$ for all metrics $g$ in the conformal
class. By construction, $\cq$ has fibre $\Bbb R_+$ and the metrics in
the conformal class are in bijective correspondence with smooth
sections of $\cq$.

Denoting by $\pi:\cq\to M$ the restriction to $\cq$ of the canonical
projection $S^2T^*M\to M$, we can view this as a principal bundle with
structure group $\Bbb R_+$. The usual convention is to rescale a metric $ g$ to
$\hat g=f^2g$. This corresponds to a  principal action given by
$\rho^s(g_x)=s^2g_x$ for $s\in\Bbb R_+$ and $g_x\in \cq_x$, the fibre
of $\cq$ over $x\in M$.

Having this, we immediately get a family of basic real line bundles
$\ce[w]\to M$ for $w\in\Bbb R$ by defining $\ce[w]$ to be the
associated bundle to $\cq$ with respect to the action of $\Bbb R_+$ on
$\Bbb R$ given by $s\cdot t:=s^{-w}t$. The usual correspondence
between sections of an associated bundle and equivariant functions on
the total space of a principal bundle then identifies the space
$\Ga(\ce[w])$ of smooth section of $\ce[w]$ with the space of all
smooth functions $f:\cq\to\Bbb R$ such that $f(\rho^s(g_x))=s^w
f(g_x)$ for all $s\in\Bbb R_+$.

Although the bundle $\ce[w]$ as we defined it depends on the choice of the 
conformal structure, it is naturally isomorphic to a density bundle (which 
is independent of the conformal structure). Recall that the bundle of
$\al$--densities is associated to the full linear frame bundle of $M$ with
respect to the $1$--dimensional representation $A\mapsto |\det(A)|^{-\al}$ of
the group $GL(n,\Bbb R)$. In particular, $1$--densities are exactly the
geometric objects that may be integrated in a coordinate--independent way on
non--orientable manifolds, while in the orientable case they may be
canonically identified with $n$--forms. To obtain the identification, recall
that any pseudo--Riemannian metric $g$ on $M$ determines a nowhere vanishing
1--density, the volume density $\vol(g)$. In local coordinates, this
density is given by $\sqrt{|\det(g_{ij})|}$, which immediately implies
that for a positive function $f$ we get $\vol(f^2g)=f^n\vol(g)$.

Consequently, any 1--density $\ph$ determines a smooth function
$\cq\to\Bbb R$ by mapping $g_x$ to $\ph(x)/\vol(g)(x)$ and obviously
this function is homogeneous of degree $-n$. This gives an
identification of the basic density bundle with $\ce[-n]$ and thus an
identification of $\ce[w]$ with the bundle of
$(-\frac{w}{n})$--densities on $M$. Hence if we have not fixed a conformal
structure in the sequel, we will switch the point of view and consider 
$\ce[w]$ as being defined as the bundle of $(-\frac{w}{n})$--densities and 
a choice of a conformal structure providing an identification of this 
density bundle with an associated bundle to $\pi:\cq\to M$.

In the sequel, we will follow the convention that adding the
expression $[w]$ to the name of any bundle indicates the tensor
product of that bundle with $\ce[w]$, so for example
$TM[-1]=TM\otimes\ce[-1]$. Clearly, sections of such weighted tensor
bundles may be viewed as equivariant sections of pullback bundles. For
example, smooth sections of $TM[w]$ are in bijective correspondence
with smooth sections $\xi$ of $\pi^*TM$ such that
$\xi(s^2g_x)=s^w\xi(g_x)$. (Recall that the fibres of $\pi^*TM$ in
$g_x$ and $s^2g_x$ may be canonically identified, so this equation
makes sense.) In particular, we may consider the tautological inclusion of
$\cq$ into $\pi^*S^2T^*M$ as a canonical section of $S^2T^*M[2]$
describing the conformal class, which gives another equivalent
description of a conformal structure.  

Of course, homogeneity along $\cq$ may as well be characterised
infinitesimally. For this, let $X$ be the fundamental vector field for
the $\Bbb R_+$--action on $\cq$,
i.e.~$X(g_x)=\frac{d}{dt}|_{t=0}\rho^{e^t}(g_x)=
\frac{d}{dt}|_{t=0}(e^{2t}g_x)$. For a function $f:\cq\to\Bbb R$ and
$w\in\Bbb R$, the equation $f(s^2g_x)=s^wf(x)$ is then clearly
equivalent to $X\cdot f=wf$. Similarly, a tensor field $t$ on $\cq$ is
called {\em homogeneous of degree\/} $w\in\Bbb R$ if and only if
$(\rho^s)^*t=s^wt$, which is equivalent to $\Cal L_Xt=wt$, where $\Cal
L$ denotes the Lie derivative. Using this, we may for example
interpret the space of smooth sections of $TM[w]$ as the quotient of
the space $\{\xi\in\frak X(\cq):[X,\xi]=w\xi\}$ by the subspace
consisting of those elements whose values in each point are
proportional to $X$.

\subsection{Tractor description of conformal structures}\label{2.2}
It is a result that goes back to E.~Cartan that conformal structures
admit a canonical normal Cartan connection, see \cite{Cartan}. More
precisely, consider $\Bbb V=\Bbb R^{n+2}$ equipped
with a non--degenerate inner product $\langle\ ,\ \rangle$ of
signature $(p+1,q+1)$. Now put $G:=O(\Bbb V)$, the orthogonal group of
$\Bbb V$, so $G\cong O(p+1,q+1)$. Furthermore, we define $P\subset G$ to be
the stabiliser of a fixed null line in $\Bbb V$. It then turns out
that $P\subset G$ is a parabolic subgroup, which may be nicely
described explicitly, see e.g.~\cite{Mike-Srni}. The relation of this
pair to conformal geometry can be described as follows: Let $\Cal C$
be the cone of nonzero null vectors in $\Bbb V$ and let $N$ be its
image in the projectivisation of $\Cal P(\Bbb V)\cong\Bbb
RP^{n+1}$. Then it is easy to see that $\langle\ ,\ \rangle$ induces a
conformal structure of signature $(p,q)$ on $N$ and $G$ acts
transitively by conformal isometries. Moreover, it turns out that this
conformal structure is flat and while $G$ does not act effectively on
$N$ (essentially since $\id$ and $-\id$ both act as the identity on
projective space) it is a two--fold covering of the group of conformal
automorphisms of $N$. Thus $N\cong G/P$ is the homogeneous flat
model of conformal geometry. 

A slight generalisation of Cartan's result may be expressed as
follows. A choice of a conformal structure on a smooth manifold $M$
gives rise to a canonical principal $P$--bundle $\Cal G\to M$ which is
endowed with a uniquely determined Cartan connection $\om\in\Om^1(\Cal
G,\frak g)$, where $\frak g=\frak o(\Bbb V)$ is the Lie algebra of $G$
and $\om$ satisfies a normalisation condition to be discussed
below. To see this in more detail (see e.g.~\cite[2.2]{luminy}) note
that the parabolic subgroup $P\subset G$ is actually related to a
grading of the Lie algebra $\frak g$ of the form $\frak
g_{-1}\oplus\frak g_0\oplus\frak g_1$. Denoting by $G_0\subset P$ the
subgroup of all elements whose adjoint action preserves this grading,
then it is elementary to verify (as outlined in \cite[2.3]{luminy})
that this group consists of $(1+n+1)\x (1+n+1)$--block matrices of the
form $\begin{pmatrix} c & 0 & 0\\ 0 & C & 0\\ 0 & 0 &
c^{-1}\end{pmatrix}$ with $0\neq c\in\Bbb R$ and $C\in O(p,q)$. The
action of such an element on $\frak g_{-1}\cong \Bbb R^n$ is given by
the standard action of $c^{-1}C$. Now one immediately verifies that 
$(c,C)\mapsto (c/|c|,c^{-1}C)$
induces an isomorphism $G_0\to\Bbb Z_2\x CO(p,q)$, where $CO(p,q)$
denotes the (pseudo--) conformal group. The inverse isomorphism is
given by $(\ep,A)\mapsto
(\ep|\det(A)|^{-1/n},\ep|\det(A)|^{-1/n}A)$. This isomorphism
intertwines the adjoint action of $G_0$ on $\frak g_{-1}$ with the
product of the trivial action of $\Bbb Z_2$ and the standard action of
$CO(p,q)$ on $\Bbb R^{p+q}$. In particular, this implies that a (first
order) $G_0$--structure is the same thing as a $CO(p,q)$--structure
and hence a conformal structure on the manifold. Now the procedure of
\cite{Cap-Schichl} applies to produce a normal Cartan connection. (Of
course in this special case there are much simpler direct
constructions of the Cartan bundle and the normal Cartan connection.)
The Cartan bundle and its normal Cartan
connection are uniquely determined by the underlying conformal
structure up to isomorphism.

While this Cartan connection is convenient from the point of view of
the equivalence problem, it is rather difficult to use it for problems
like finding invariants of conformal structures or conformally
invariant differential operators. To deal with such problems, it is
often more efficient to switch to the description of conformal
structures via the so--called standard tractor bundle and its
canonical linear connection. Starting from the Cartan bundle and the
Cartan connection, the standard tractor bundle $\Cal T\to M$ is simply
the associated bundle $\Cal G\x_{P}\Bbb V$. By construction, this
bundle carries a canonical metric of signature $(p+1,q+1)$. The
distinguished null line in $\Bbb V$ used to define $P$ leads to a
subbundle $\Cal T^1\subset\Cal T$ whose fibres are null lines and
which is easily seen to be isomorphic to $\ce[-1]$. Furthermore, it
turns out that the Cartan connection $\om$ induces a linear connection
$\nabla$ on $\Cal T$, the so--called normal standard tractor
connection, see \cite{tractors}. Having these data at hand, one may
then compute the fundamental D--operator on $\Cal T$, see
\cite[section 3]{tractors}, which in turn leads to the so called
tractor D--operator, see \cite{gosrni} or \cite[section
3]{luminy}. These operators have been successfully applied to the
construction of conformally invariant differential operators,
conformal invariants and other topics in conformal geometry, see
e.g.~\cite{BEG,nonlocal,conf-invar}. In summary, having an explicit
knowledge of the standard tractor bundle, the tractor metric and the
normal standard tractor connection, one immediately gets a large
number of tools for dealing with problems in conformal geometry.

It is an idea going back to the work of T.~Thomas in the 20's to use
the standard tractor bundle and its canonical connection as an
alternative approach to conformal geometry. The precise relation
between these data and the Cartan bundle and Cartan geometry was
completely clarified (in a much more general setting) in
\cite{tractors}. Specialised to conformal standard tractors, this goes
as follows: Suppose that $M$ is a smooth manifold, and that $\Cal T\to
M$ is a real rank $n+2$ vector bundle endowed with a bundle metric $h$
of signature $(p+1,q+1)$ and an injective bundle map $\ce[-1]\to\Cal
T$, whose image $\Cal T^1$ is null with respect to $h$. Suppose
further that the $ \Cal T$ admits a tractor connection $\nabla$ in the
sense of \cite[2.5]{tractors}. By definition, this means that $\nabla$
is a non--degenerate $\frak o(\Bbb V)$--connection. The condition that
$\nabla$ is a $\frak o(\Bbb V)$ connection is easily seen to be
equivalent to $\nabla$ preserving the bundle metric $h$. On the other
hand, $\Cal T$ has a filtration of the form $\Cal T\supset \Cal
T^0\supset \Cal T^1$, where $\Cal T^0:=(\Cal T^1)^\perp$. This
immediately implies that for connections preserving $h$ the
non--degeneracy condition from \cite[2.5]{tractors} is equivalent to
the condition that for any $x\in M$ and $\xi\in T_xM$ there is a
smooth section $\si$ of $\Cal T^1$ such that
$\nabla_\xi\si(x)\notin\Cal T^1_x$.

Given the data $(\Cal T,h,\nabla)$ as above, we can now recover an
underlying conformal structure of signature $(p,q)$ on $M$: First, let
$\si_0$ be a locally non--vanishing section of $\Cal T^1$. Then
$h(\si_0,\si_0)=0$ and thus $0=\xi\cdot
h(\si_0,\si_0)=2h(\nabla_\xi\si_0,\si_0)$ for all $\xi\in\frak
X(M)$. This immediately implies that for any smooth function $f$ we
get $h(\nabla_\xi(f\si_0),\si_0)=0$, and since locally any smooth
section of $\Cal T^1$ can be written in the form $f\si_0$, we conclude
that $\nabla_\xi\si\in\Ga(\Cal T^0)$ for all $\xi\in\frak
X(M)$ and $\si\in\Ga(\Cal T^1)$. Now consider the map which maps
$(\xi,\si)$ to the class of $\nabla_\xi\si$ in $\Ga(\Cal
T^0/\Cal T^1)$. This is obviously bilinear over smooth
functions, and thus induced by a bundle map $TM\otimes\ce[-1]\to\Cal
T^0/\Cal T^1$, which by the non--degeneracy assumption is
injective on each fibre, so since both bundles have rank $n$, we
obtain a bundle isomorphism $\Cal T^0/\Cal T^1\cong
TM[-1]$. On the other hand, since the restriction of $h$ to $\Cal
T^0$ is degenerate with null space $\Cal T^1$, $h$ induces a
non--degenerate bundle metric of signature $(p,q)$ on $\Cal
T^0/\Cal T^1$, and thus gives rise to a section of
$S^2T^*M[2]$, i.e.~a conformal structure on $M$. We say $(\Cal
T,h,\nabla)$ is a standard tractor bundle corresponding to this
conformal structure. Conversely beginning with a conformal structure on
$ M$ there are ways (see e.g.  \cite{BEG,luminy}) to directly
construct standard tractor bundles for the given conformal structure.

Next, one may recover the Cartan bundle from the standard tractor
bundle: For $x\in M$ define $\Cal G_x$ to be the set of all orthogonal
maps $\Bbb V\to\Cal T_x$ which in addition map the distinguished null
line to $\Cal T^1_x$. By assumption, such maps exist, and composition
from the right defines a transitive free right action of $P$ on $\Cal
G_x$. Now the union $\Cal G:=\cup_{x\in M}\Cal G_x$ may be naturally
viewed as a subbundle of the frame bundle of $\Cal T$, whence it
obtains its smooth structure and the $P$--action from above makes it
into a $P$--principal bundle. By construction, we have $\Cal T=\Cal
G\x_P\Bbb V$ and the metric $h$ and the subbundle $\Cal T^1$ are
obtained by carrying over the respective data from $\Bbb V$. In the
language of \cite[section 2]{tractors}, this means that $\Cal T\to M$
is a standard tractor bundle and $\Cal G$ is an adapted frame bundle
for $\Cal T$. (The adjoint tractor bundle lurking in the background is
the bundle $\frak o(\Cal T)$ of skew symmetric endomorphisms of $\Cal
T$).

Now by \cite[Theorem 2.7]{tractors} there is a bijective
correspondence between tractor connections on $\Cal T$ and Cartan
connections on the adapted frame bundle $\Cal G$. To recognise the
normal tractor connection among all tractor connections, one notes
that by \cite[Proposition 2.9]{tractors}, the curvature $R$ of a
tractor connection $\nabla$ is given by the action of the curvature of
the corresponding Cartan connection. Now in the special case of
conformal structures, the general Lie theoretic normalisation
condition on Cartan connections used in \cite{tractors} can be
simplified considerably. First of all, any normal Cartan connection in
the conformal case is torsion free, which simply means that the action
of $R(\xi,\eta)$ on the standard tractor bundle preserves the
subbundle $\Cal T^1\subset\Cal T$ for all $\xi$ and $\eta$. If this
condition is satisfied, then $R(\xi,\eta)$ induces an endomorphism
$W(\xi,\eta)$ of $\Cal T^0/\Cal T^1\cong TM[-1]$, so we may as well
view $W$ as a section of $\La^2T^*M\otimes L(TM,TM)$. Using this and
taking into account the description of $\partial^*$ in the proof of
\cite[Proposition 4.3]{tractors} and the formula for the algebraic
bracket in the conformal case in \cite[2.3]{luminy} one concludes that
the normalisation condition on the curvature of a standard tractor
connection is equivalent to vanishing of the Ricci--type contraction
of $W$. Of course, uniqueness of the Cartan bundle and Cartan
connection implies that $\Cal T$ together with the subbundle $\Cal
T^1$, the metric $h$ and the normal standard tractor connection
$\nabla$ is uniquely determined by the underlying conformal structure
up to isomorphism. Summarising, we obtain
\begin{thm*}
  (1) Let $M$ be a smooth manifold of dimension $n\geq 3$. Suppose
  that $\Cal T\to M$ is a rank $n+2$ real vector bundle endowed with a
  bundle metric $h$ of signature $(p+1,q+1)$, an injective bundle map
  $\ce[-1]\to\Cal T$ with null image $\Cal T^1\subset\Cal T$ and a
  linear connection $\nabla$ such that $\nabla h=0$ and, for any $x\in
  M$ and any $\xi\in T_xM$, there is a smooth section $\si\in\Ga(\Cal
  T^1)$ such that $\nabla_\xi\si(x)\notin\Cal T^1_x$. 
 Then $(\Cal T^1)^\perp/\Cal T^1\cong
  TM[-1]$, and $(\Cal T,h,\nabla)$ is a standard tractor bundle for
  the conformal structure defined by the restriction of $h$ to $(\Cal
  T^1)^\perp/\Cal T^1\x (\Cal T^1)^\perp/\Cal T^1$.

\noindent
(2) The tractor connection $\nabla$ on $\Cal T$  is normal if
and only if its curvature $R$ has the properties that $R(\xi,\eta)(\Cal
T^1)\subset\Cal T^1$ and the Ricci--type contraction of the element
$W\in\Ga(\La^2T^*M\otimes L(TM,TM))$, as described above, vanishes.
If this is the case, then $(\Cal T,\Cal T^1,h,\nabla)$ is uniquely
determined by the underlying conformal structure up to isomorphism. 
\end{thm*}

\subsection{Ambient manifolds and ambient metrics}\label{2.3}
In \cite{F-G}, Ch.~Fefferman and C.R.~Graham have initiated a project
to study conformal structures using the so--called ambient metric
construction. The idea of that construction is to mimic the flat
metric on the vector space $\Bbb V$ in the case of the homogeneous
model as described in \ref{2.2} above. The null cone $\Cal C$ may be
viewed as the image of the inclusion $\ce[-2]\to S^2T^*M$ provided by
the conformally flat structure, so half of this null cone may be
identified with the bundle $\cq$ of metrics in the conformal
class. Now one starts with the bundle $\cq\to M$ of metrics defining
an arbitrary conformal structure. Then in \cite{F-G} it is shown that
there is a Riemannian metric of signature $(p+1,q+1)$ on $\cq\x
(-1,1)$ (defined locally around $\cq$) whose Ricci curvature vanishes
to a certain order (depending on the dimension) along $\cq$. Moreover,
the corresponding jet of this metric along $\cq$ is unique in a
certain sense. This construction has been applied in \cite{GJMS} to
prove the existence of certain conformally invariant powers of the
Laplacian. Using $\cq\x (-1,1)$ is slightly misleading, one could equally
consider a germ along $\cq$ of an (unspecified) ambient
manifold endowed with a free $\Bbb R_+$--action.  Moreover, as we
shall see later on, a much weaker normalisation condition on an ambient
metric than the one used in \cite{F-G} is sufficient to get the
relation to standard tractors. Thus, we will start our discussion with
a general version of ambient manifolds and ambient metrics.

Note first, that on any manifold endowed with a free action of $\Bbb
R_+$, one has the notion of homogeneity of tensor fields as described
in \ref{2.2} above, which can equivalently be characterised
infinitesimally. 

\subsection*{Definition}
Let $\pi:\cq\to M$ be a conformal structure. An {\em ambient manifold\/}
is a smooth $(n+2)$-manifold $\aM$ endowed with a free $\Bbb R_+$--action
$\rho$ on $\aM$ and a $\Bbb R_+$--equivariant embedding
$\iota:\cq\to\aM$. 

If $\iota:\cq\to\aM$ is an ambient manifold, then an {\em ambient
metric\/} is a pseudo--Riemannian metric $\h$ of signature $(p+1,q+1)$
on $\aM$ such that the following condition hold:

\noindent
(i) The metric $\h$ is homogeneous of degree 2 with respect to the
$\Bbb R_+$--action, i.e.~if $\X\in\frak X(\aM)$ denotes the
fundamental field generating the $\Bbb R_+$--action and $\Cal L_{\sX}$
denotes the Lie derivative by $\X$, then we have $\Cal L_{\sX}\h=2\h$.

\noindent
(ii) For $u=g_x\in \cq\subset\aM$ and $\xi,\eta\in T_u\cq$, we have
$\h(\xi,\eta)=g_x(T\pi\cdot\xi,T\pi\cdot\eta)$. 

\smallskip

Since the action of $\Bbb R_+$ on $\aM$ extends the action on
$\cq$, we will denote both actions by the symbol $\rho$ and we use $\X$
to denote the fundamental field for both actions. Moreover, we will
usually view $\cq$ as a submanifold of $\aM$ and suppress the
embedding $\iota$.

Since we will frequently have to deal with the question of vanishing
of tensor fields along $\cq$ to some order, we collect some
information on that. A tensor field $t$ on $\aM$ is said to vanish
along $Q$ to order $\ell\geq 1$ if and only if $t|_Q=0$ and for any
integer $k<\ell$ and arbitrary vector fields
$\xi_1,\dots,\xi_k\in\frak X(\tilde M)$ the iterated Lie derivative
$\Cal L_{\xi_k}\cdots\Cal L_{\xi_1}t$ vanishes along
$\cq$. Equivalently, one may require all iterated covariant
derivatives $\annd_{\xi_k}\cdots\annd_{\xi_1}t$ to vanish along
$\cq$. The tensor field $t$ is said to vanish to infinite order along
$\cq$ if it vanishes to order $\ell$ for all $\ell\in\Bbb N$. If we
choose any defining function $r$ for $\cq$, i.e.~a smooth real valued
function defined locally around $\cq$ such that $\cq=r^{-1}(0)$ and
$dr$ does not vanish in any point of $\cq$, then any tensor field $t$
that vanishes along $\cq$ may be written as $t=rt'$ for some tensor
field $t'$ of the same type as $t$. Inductively, one sees that $t$
vanishes to order $\ell$ along $\cq$ if and only if $t=r^\ell t'$ for
some tensor field $t'$. Thus, we will use the notation $t=O(r^\ell)$
to indicate that $t$ vanishes to order $\ell$ along $Q$. 

There are some points we should make that are particular to the 
 case of sections of $\otimes^sT^*M$. This is especially relevant for 
the case of differential forms. On the one
hand, in this case vanishing to order $\ell$ along $\cq$ can be
equivalently expressed as vanishing of $\xi_k\cdots\xi_1\cdot
t(\eta_1,\dots,\eta_s)$ for arbitrary $k<\ell$ and vector fields
$\xi_i$ and $\eta_j$. On the other hand, for tensor fields $t$ of that
type there is the weaker condition that {\em tangential components of
$t$\/} vanish to some order along $\cq$. One says that the tangential
components of $t$ vanish along $\cq$ if $\iota^*t=0$, where
$\iota:\cq\to\aM$ is the inclusion. Equivalently,
$t(\eta_1,\dots,\eta_s)$ has to vanish along $\cq$ if for any $u\in Q$
and any $j$ one has $\eta_j(u)\in T_u\cq\subset T_u\aM$. We say that
the tangential components of $t$ vanish to order $\ell$ if and only if
for $k<\ell$, arbitrary vector fields $\xi_1,\dots,\xi_k\in\frak
X(\aM)$ and vector fields $\eta_1,\dots,\eta_s\in\frak X(\aM)$ such
that each $\eta_j|_{\cq}$ is tangent to $\cq$, the function
$\xi_k\cdots\xi_1\cdot t(\eta_1,\dots,\eta_s)$ vanishes along
$\cq$. In this case, however, one may not replace this by a condition
on Lie derivatives or covariant derivatives, since a Lie derivative or
covariant derivative of a vector field whose restriction to $\cq$ is
tangent to $\cq$ in general does not have the same property.

The normalisation conditions on ambient metrics used in \cite{F-G} are
based on the Ricci curvature of the ambient metric $\h$. However, to
get the relation to standard tractors, we need a different
condition. We shall show in \ref{2.6} that this condition is a
consequence of the weakest possible condition on the Ricci
curvature. The condition we need is based on the one--form $\al$ dual
to the generator $\X$ of the $\Bbb R_+$--action,
i.e.~$\al(\xi)=\h(\X,\xi)$. Notice that since $T\pi\cdot\X=0$,
condition (ii) in the definition of an ambient metric implies that
$\iota^*\al=0$. Thus, we also have $0=d\iota^*\al=\iota^*d\al$, so the
tangential components of $d\al$ vanish along $\cq$.

Expanding the homogeneity condition $\Cal L_{\sX}\h=2\h$, we get
${\X}\cdot\h(\xi,\eta)-\h([{\X},\xi],\eta)-\h(\xi,[{\X},\eta])
=2\h(\xi,\eta)$, and rewriting the Lie brackets in terms of covariant
derivatives, we obtain 
$$
\h(\annd_\xi\X,\eta)+\h(\xi,\annd_\eta\X)=2\h(\xi,\eta),
$$
which says that $\X$ is a conformal Killing field of dilation type. On
the other hand, by definition of the exterior derivative, we get
$d\al(\xi,\eta)=\xi\cdot\h(\X,\eta)-\eta\cdot\h(\X,\xi)-\h(\X,[\xi,\eta])$,
and expanding the right hand side of this in terms of covariant
derivatives gives
$$
\h(\annd_\xi\X,\eta)-\h(\xi,\annd_\eta\X)=d\al(\xi,\eta),
$$
which just expresses the fact that the exterior derivative of a
differential form is obtained by alternating the covariant
derivative. Putting these two equations together, we obtain 
\begin{equation}\label{nX}
\h(\annd_\xi\X,\eta)=\h(\xi,\eta)+\frac{1}{2}d\al(\xi,\eta). 
\end{equation}
This equation shows that $d\al=O(r^\ell)$ implies
$\annd_\xi\X=\xi+O(r^\ell)$ for any vector field $\xi$. On the other
hand, if $\xi$ is homogeneous of degree $w$, then
$w\xi=[\X,\xi]=\annd_{\sX}\xi-\annd_\xi\X$, which shows that
$\annd_{\sX}\xi=(w+1)\xi+O(r^\ell)$ provided that
$d\al=O(r^\ell)$. Suppose next that $i_{\sX}d\al|_{\cq}=0$, where $i$
denotes the insertion operator. Since $\al$ is obviously homogeneous
of degree two, this implies that $2\al|_{\cq}=\Cal
L_{\sX}\al|_{\cq}=di_{\sX}\al|_{\cq}$. Thus we see that assuming
$i_{\sX}d\al|_{\cq}=0$, the function $r:=\frac{1}{2}\h(\X,\X)$
satisfies $dr|_{\cq}=\al|_{\cq}$, so we get a canonical defining
function in this case.

\subsection{The standard tractor bundle and connection induced by an
ambient metric}\label{2.4}
Let $\pi:\cq\to M$ be a conformal structure, $\aM\supset \cq$ an
ambient manifold and $\h$ and ambient metric on $\aM$. Throughout this
subsection we assume that $\h$ has the property that the one--form
$\al$ dual to $\X$ satisfies $d\al|_Q=0$.

Consider the restriction $T\aM|_\cq$ of the ambient tangent bundle
to $\cq$ and define an action of $\Bbb R_+$ on this space by $s\cdot
\xi:=s^{-1}T\rho^s\cdot\xi$. This is compatible with the $\Bbb R_+$ 
action on $\cq$, so defining $\Cal T$ to be the
quotient $(T\aM|_\cq)/\Bbb R_+$, we immediately see that this is a
vector bundle over $\cq/\Bbb R_+=M$, and the fibre dimension of
this bundle is $n+2$. Moreover, by construction, there is a bijective
correspondence between the space $\Ga(\Cal T)$ of smooth sections of
$p:\Cal T\to M$ and the space of ambient vector fields along $\cq$
(i.e.~sections of $T\aM|_\cq\to \cq$) which are homogeneous of degree
$-1$, or equivalently satisfy $[\X,\xi]=-\xi$. 

The fact that the ambient metric $\h$ is homogeneous of degree $2$
immediately implies that for vector fields $\xi$ and $\eta$ on $\aM$
which are homogeneous of degree $w$ and $w'$, respectively, the
function $\h(\xi,\eta)$ is homogeneous of degree $w+w'+2$. In
particular, applying $\h$ to the vector fields corresponding to two
sections of $\Cal T$, the resulting function on $\cq$ is constant on
$\Bbb R_+$ orbits, and thus descends to a smooth function on
$M$. Hence $\h$ descends to a smooth bundle metric $h$ of signature
$(p+1,q+1)$ on $\Cal T$. 

The bundle metric $h$ immediately leads to a filtration of the bundle
$\Cal T$: Since the vertical tangent bundle of $\pi:\cq\to M$ is
stable under the $\Bbb R_+$--action, it gives rise to a distinguished
line bundle $\Cal T^1\subset\Cal T$. By construction, sections of this
subbundle correspond to ambient vector fields along $\cq$, which are
of the form $f{\X}$ for some smooth function $f:\cq\to\Bbb R$, and in
order that this is a section of $\Cal T$, the function $f$ must be
homogeneous of degree $-1$. Thus, mapping $f$ to $f{\X}$ defines an
isomorphism $\ce[-1]\cong\Cal T^1$. On the other hand, we have already
observed that $\h({\X},{\X})=0$ along $\cq$. Hence, defining $\Cal
T^0$ to be the orthogonal complement of $\Cal T^1$ with respect to
$\h$, we see that $\Cal T^0\subset\Cal T$ is a smooth subbundle of
rank $n+1$ and $\Cal T^1\subset\Cal T^0$. To identify the quotient
$\Cal T/\Cal T^0$, we observe that for any section $s\in\Ga (\Cal T)$
with corresponding vector field $\xi$ along $\cq$, we get a function
$\h(\xi,{\X})$, which is homogeneous of degree one. By construction,
this vanishes if and only if $s$ has values in $\Cal T^0$, so it
induces an isomorphism $\Cal T/\Cal T^0\cong\ce[1]$ of vector bundles.

Finally, assume that $\xi$ is a vector field on $M$ and
$f\in\Ga(\ce[-1])$ is a smooth section, i.e.~a function $\cq\to\Bbb R$
homogeneous of degree $-1$. Then we may lift $\xi$ to an ambient
vector field $\tilde\xi$ along $\cq$, which is homogeneous of degree
zero, and this lift is unique up to adding fields of the form $\ph
{\X}$ with $\ph:\cq\to\Bbb R$ homogeneous of degree zero. Then
$f\tilde\xi$ is a section of $\Cal T$ and by property (ii) of $ \h$ we
have $\h(f\tilde\xi,{\X})=0$, whence $f\tilde\xi\in\Ga(\Cal
T^0)$. Moreover, the class of $f\tilde\xi$ in $\Cal T^0/\Cal T^1$ is
independent of the choice of the lift $\tilde\xi$. Hence we obtain a
bundle map $TM[-1]\to\Cal T^0/\Cal T^1$, which is obviously injective
in each fibre, so since both bundles have the same rank, we conclude
that $\Cal T^0/\Cal T^1\cong TM[-1]$, which implies that $p:\Cal T\to
M$ is a candidate for a conformal standard tractor bundle. Notice that
up to now, we have only used properties (i) and (ii) of the ambient
metric.

Next, let $\annd$ be the Levi--Civita connection of $\h$. The
fact that $\annd$ is torsion free and $\annd\h=0$ imply
the well know global formula
$$
\begin{aligned}
2\h(\annd_\xi\eta,\zeta)=
&\xi\cdot\h(\eta,\zeta)+\eta\cdot\h(\xi,\zeta)-\zeta\cdot\h(\xi,\eta)\\ 
+&\h([\xi,\eta],\zeta)-\h([\xi,\zeta],\eta)-\h([\eta,\zeta],\xi) 
\end{aligned}
$$
for all vector fields $\xi,\eta,\zeta\in\frak X(\aM)$.  Observe that
if $\xi\in\frak X(\aM)$ is homogeneous of degree $w$ and $f:\aM\to\Bbb
R$ is homogeneous of degree $w'$, then the equation ${\X}\cdot
\xi\cdot f=[\X,\xi]\cdot f+\xi\cdot {\X}\cdot f$ shows that the
function $\xi\cdot f$ is homogeneous of degree $w+w'$.  Hence choosing
the three vector fields in the above formula to be homogeneous of
degrees $w$, $w'$, and $w''$, we see that any summand on the right
hand side is homogeneous of degree $w+w'+w''+2$, which immediately
implies that $\annd_\xi\eta$ is homogeneous of degree $w+w'$. In
particular, if $\xi$ is invariant, i.e.~homogeneous of degree zero,
then $\annd_\xi\eta$ has the same homogeneity as $\eta$.

On the other hand, we have already observed in the end of \ref{2.3}
above, that $d\al|_Q=0$, i.e.~$d\al=O(r)$ implies that
$\annd_\xi\X=\xi+O(r)$, and $\annd_{\sX}\xi=(w+1)\xi+O(r)$ if $\xi$ is
homogeneous of degree $w$. In particular, $\annd_{\sX}\xi|_Q=0$ for
$\xi$ homogeneous of degree $-1$. Using this, we can now show that
$\annd$ descends to a linear connection $\nabla$ on the bundle $\Cal
T$. Suppose that $s\in\Ga(\Cal T)$ is a section corresponding to the
ambient vector field $\tilde\eta$ along $\cq$ and that $\xi\in\frak
X(M)$ is a vector field. As before, we may lift $\xi$ to an ambient
vector field $\tilde\xi$ along $\cq$, which is unique up to adding
terms of the form $\ph {\X}$ with $\ph$ homogeneous of degree
zero. Extend $\tilde\eta $ to a homogeneous field on $ \aM$ and
observe that since the flow lines of $\tilde\xi$ are contained in
$\cq$, it follows that, along $ \cq$, the ambient vector field
$\annd_{\tilde\xi}\tilde\eta$ is independent of the extension of $
\tilde\eta$. Since $\annd_{\sX}\tilde\eta|_Q=0$, we conclude that
$\annd_{\tilde\xi}\tilde\eta$ depends only on $\xi$ and not on the
lift $\tilde\xi$. Moreover, from above we know that
$\annd_{\tilde\xi}\tilde\eta$ is homogeneous of degree $-1$, so it
corresponds to a section of $\Cal T$, which we denote by $\nabla_\xi
s$. One immediately verifies that this defines a linear connection
$\nabla$ on $\Cal T$, which by construction is compatible with the
bundle metric $h$.

To verify that $\nabla$ is a tractor connection on $\Cal T$, we thus
only have to verify the non--degeneracy condition from \ref{2.2},
which is very easy: For a section $s\in\Ga(\Cal T^1)$, the
corresponding ambient field is of the form $f{\X}$ with $f:\cq\to\Bbb
R$ homogeneous of degree $-1$. For a lift $\tilde\xi$ of a vector
field $\xi\in\frak X(M)$ as above, we get
$\annd_{\tilde\xi}f{\X}=(\tilde\xi\cdot f){\X}+f\annd_{\tilde\xi}{\X}$
and the second summand equals $f\tilde\xi$ along $\cq$. In particular,
we see that $\nabla_\xi s\in\Ga(\Cal T^0)$, and the image of this
section in $\Ga(\Cal T^0/\Cal T^1)$ is simply the element $f\xi$,
which implies that $\nabla$ is a tractor connection on $\Cal T$. Thus
we have proved:

\begin{prop*}
  Let $\pi:\cq\to M$ be a conformal structure on a smooth manifold
  $M$, $\aM$ an ambient manifold and $\h$ an ambient metric, and let
  $\al\in\Om^1(\aM)$ be the one--form dual to the infinitesimal
  generator of the $\Bbb R_+$--action on $\aM$. Then for the $\Bbb
  R_+$ action on $T\aM|_\cq$ defined above, $\h$ descends to a bundle
  metric $h$ on $\Cal T:=(T\aM|_\cq)/\Bbb R_+$. If $d\al|_{\cq}=0$,
  then the Levi--Civita connection of $\h$ descends to a tractor
  connection on $\Cal T$ which preserves $ h$. This together with the
  filtration induced by the vertical subbundle means $\Cal T$ is a
  conformal standard tractor bundle.
\end{prop*}

\subsection{The normalisation condition}\label{2.5}  
Let us assume that $\h$ is an ambient metric on an ambient manifold
$\aM$ for a given conformal structure on $M$ such that the one--form
$\al$ dual to the infinitesimal generator $\X$ of the $\Bbb
R_+$--action on $\aM$ has the property that $d\al|_{\cq}=0$. Then we
have the induced conformal standard tractor bundle $(\Cal
T,h,\nabla)$. Now it is almost obvious that the curvature $\aR$ of
$\h$ descends to curvature $R$ of the tractor connection $\nabla$.
Indeed, choosing invariant lifts $\tilde\xi$ and $\tilde\eta$ for
vector fields $\xi,\eta\in\frak X(M)$, and considering the
(homogeneous of degree $-1$) ambient vector field $\zeta$ along $\cq$
corresponding to a section $s\in\Ga(\Cal T)$, the ambient vector field
$\annd_\xi\annd_\eta\zeta$ corresponds to the section
$\nabla_\xi\nabla_\eta s\in\Ga(\Cal T)$. Moreover,
$[\tilde\xi,\tilde\eta]$ is an invariant lift of $[\xi,\eta]$, which
immediately implies that $\aR(\tilde\xi,\tilde\eta)\zeta$ corresponds
to the section $R(\xi,\eta)s$ of $\Cal T$.

Hence to understand the tractor curvature $R$, we have to analyse the
ambient curvature $\aR$. For later use, we work in a more general
setting and prove more specific results than are required for the
verification of the normalisation condition.

\begin{prop*}
Let $\h$ be an ambient metric on $\aM$, $\aR$ the Riemann curvature of
$\h$, $\X\in\frak X(\aM)$ the infinitesimal generator of the $\Bbb
R_+$--action on $\aM$ and $\al\in\Om^1(\aM)$ the one--form dual to
$\X$. Then we have:
$$\h(\aR(\xi,\eta)\X,\zeta)=-\h(\aR(\xi,\eta)\zeta,\X)=
    \h(\aR(\X,\zeta)\xi,\eta)=-\tfrac{1}{2}(\annd d\al)(\zeta,\xi,\eta),
$$
for all vector fields $\xi,\eta,\zeta\in\frak X(\aM)$. In particular,
if $d\al=O(r^\ell)$ for some $\ell\geq 1$, then this expression is
$O(r^{\ell-1})$ and it is $O(r^\ell)$ if either $\zeta$ or both $\xi$
and $\eta$ have the property that the restriction to $\cq$ is
tangent to $\cq$. Hence, if $d\al=O(r^\ell)$, then tangential
components of $d\al$ vanish to order $\ell+1$ along $\cq$.
\end{prop*}
\begin{proof}
The equality of the first three expressions follows from standard
symmetries of the curvature of a pseudo--Riemannian metric. Now we
compute 
$$
\h(\annd_\xi\annd_\eta\X,\zeta)=\xi\cdot\h(\annd_\eta\X,\zeta)-
\h(\annd_\eta\X,\annd_\xi\zeta), 
$$
and inserting formula \eqref{nX} from \ref{2.3} in both summands, we
obtain
$$
\h(\annd_\xi\eta,\zeta)+\tfrac{1}{2}\xi\cdot
d\al(\eta,\zeta)-\tfrac{1}{2}d\al(\eta,\annd_\xi\zeta). 
$$
Taking the alternation of this in $\xi$ and $\eta$ and subtracting
$\h(\annd_{[\xi,\eta]}\X,\zeta)=\h([\xi,\eta],\zeta)+
\tfrac{1}{2}d\al([\xi,\eta],\zeta)$ we may expand the Lie bracket into
covariant derivatives which implies that all terms involving $\h$
cancel, and we are left with the expression 
$$
\tfrac{1}{2}\bigg(\xi\cdot d\al(\eta,\zeta)-
\eta\cdot d\al(\xi,\zeta)-d\al([\xi,\eta],\zeta)-
d\al(\eta,\annd_\xi\zeta)+d\al(\xi,\annd_\eta\zeta)
\bigg)
$$
for $\h(\aR(\xi,\eta)\X,\zeta)$.  Expanding the equation
$0=d(d\al)(\xi,\eta,\zeta)$ we may rewrite $\xi\cdot
d\al(\eta,\zeta)-\eta\cdot d\al(\xi,\zeta)- d\al([\xi,\eta],\zeta)$ as
$-\zeta\cdot
d\al(\xi,\eta)-d\al([\xi,\zeta],\eta)+d\al([\eta,\zeta],\xi)$, and
expressing the Lie brackets as commutators of covariant derivatives,
we arrive at the claimed formula for $\h(\aR(\xi,\eta)\X,\zeta)$.

By definition, $(\annd d\al)(\zeta,\xi,\eta)=\zeta\cdot
d\al(\xi,\eta)-d\al(\annd_\zeta\xi,\eta)-d\al(\xi,\annd_\zeta\eta)$,
and if $d\al=O(r^\ell)$, the the last two terms visibly are
$O(r^\ell)$, while the first is $O(r^{\ell-1})$. If in addition
$\zeta|_{\cq}$ is tangent to $\cq$, the equation
$\zeta'\cdot\zeta\cdot d\al(\xi,\eta)=\zeta\cdot\zeta'\cdot
d\al(\xi,\eta)+[\zeta',\zeta]\cdot d\al(\xi,\eta)$ shows that the
first summand is $O(r^\ell)$, too. On the other hand, if both
$\xi|_{\cq}$ and $\eta|_{\cq}$ are tangent to $\cq$, then by the
Bianchi identity, we get
$\h(R(\xi,\eta)\zeta,\X)=-\h(R(\zeta,\xi)\eta,\X)-\h(R(\eta,\zeta)\xi,\X)$,
and from above we know that both terms of the right hand side are
$O(r^\ell)$. Turning around the argument, we see now that if
$\xi|_{\cq}$ and $\eta|_{\cq}$ are tangent to $\cq$, then for any
vector field $\zeta$, the function $\zeta\cdot d\al(\xi,\eta)$ is
$O(r^\ell)$, whence tangential components of $d\al$ vanish to order
$\ell+1$ along $\cq$.  
\end{proof}

Using these results we can now prove:

\begin{thm*}
Let $\h$ be an ambient metric such that the one--form $\al$ dual to
the infinitesimal generator of the $\Bbb R_+$--action satisfies
$d\al|_{\cq}=0$. Then the standard tractor bundle $(\Cal T,h,\nabla)$
induced by $\h$ is normal if and only if the tangential
components of the Ricci curvature $\Ric(\h)$ vanish along $\cq$.
\end{thm*}
\begin{proof}
From the above Proposition we see that $d\al=O(r)$ implies that if
$\tilde\xi|_{\cq}$ and $\tilde\eta|_{\cq}$ are tangent to $\cq$, then
$\h(\aR(\tilde\xi,\tilde\eta)\X,\zeta)$ and
$\h(\aR(\tilde\xi,\tilde\eta)\zeta,\X)$ vanish along $\cq$. Applied to
invariant lifts of vector fields $\xi,\eta\in\frak X(M)$, the first
equation exactly means that $R(\xi,\eta)$ vanishes on $\Cal T^1$,
while the second equation says that $R(\xi,\eta)$ maps $\Cal T$ to
$\Cal T^0$. In particular, $R(\xi,\eta):\Cal T\to\Cal T$ is filtration
preserving for all $\xi$, $\eta$.

Hence $\aR(\xi,\eta)$ induces an endomorphism of $\Cal
T^0/\Cal T^1$, which we may as well view as an endomorphism of
$TM$. From \ref{2.2} we know that normality of the tractor connection
is equivalent to vanishing of the Ricci--type contraction of the
resulting operator $W\in\Ga(\La^2T^*M\otimes L(TM,TM))$. The value of
this contraction at a point $x\in M$ on tangent vectors $\xi,\eta\in
T_xM$ can be computed as $\sum_{i=1}^n\ph_i(W(\xi_i,\xi)\eta)$, for a
basis $\{\xi_1,\dots,\xi_n\}$ of $T_xM$ with dual basis
$\{\ph_1,\dots,\ph_n\}$ of $T^*_xM$.

Choosing a point $u$ in $\cq$ over $x$ and lifts $\tilde\xi_i$,
$\tilde\xi$ and $\tilde\eta$ of the vector fields involved, we may
compute $\ph_i$ as $\h(\_,\tilde\eta_i)$, where the tangent
vectors $\tilde\eta_1,\dots,\tilde\eta_n\in T_u\cq\subset T_u\aM$
are defined by $\h(\tilde\eta_i,\tilde\xi_j)=\delta_{ij}$. (Note that
$\tilde\eta_i\in T_u\cq$ implies $\h(\tilde\eta_i,\X)=0$). Hence our
contraction applied to $\xi$ and $\eta$ corresponds to 
$$
\sum_{i=1}^n\h(\aR(\tilde\xi_i,\tilde\xi)\tilde\eta,\eta_i).
$$ Now let $Y\in T_u\cq$ be the unique null tangent vector such that
$\bar h(\X,Y)=1$ and $\h(\tilde\xi_i,Y)=0$ for all $i$. Then clearly
$\{\X,\tilde\xi_1,\dots,\tilde\xi_n,Y\}$ is a basis of $T_u\aM$ with
dual basis (with respect to $\h$) given by
$\{Y,\tilde\eta_1,\dots,\tilde\eta_n,\X\}$.  But from the above
proposition, we know that $\h(\aR(\X,\tilde\xi)\tilde\eta,Y)$ and
$\h(\aR(Y,\tilde\xi)\tilde\eta,\X)$ vanish along $\cq$, since the
restrictions of $\tilde\xi$ and $\tilde\eta$ to $\cq$ are tangent to
$\cq$. Adding these two summands to the above sum, we by definition
get $\Ric(\h)(\tilde\xi,\tilde\eta)(u)$, where $\Ric$ denotes the
Ricci curvature of $\h$, which implies the result.
\end{proof} 

Note that it follows immediately from the theorem that the normal
tractor curvature $R$ is induced by  the
curvature of any ambient metric $\h$ which has the property that
$d\al$ and tangential components of $\Ric(\h)$ vanish along $\cq$.

\subsection{}\label{2.6}
We next want to show that an ambient metric $\h$, such that the
tangential components of $\Ric(\h)$ vanish along $\cq$, automatically
satisfies $d\al|_{\cq}=0$, where $\al$ is the one--form dual to the
infinitesimal generator $\X$ of the $\Bbb R_+$--action. In particular,
this implies that any ambient metric satisfying the Ricci condition
can be used to construct the normal standard tractor bundle. 

Let us again start with an arbitrary ambient metric $\h$ on an ambient
manifold $\aM$ for $\pi:\cq\to M$. Let us first choose appropriate
dual frames defined locally around a point in $\cq$. Note that along
$\cq$, the tangent spaces of $\cq$ are orthogonal to $\X$. Thus,
locally around a point $u_0\in\cq$, we may choose ambient vector
fields $\xi_i\in\frak X(\aM)$ for $i=1,\dots, n$, which are
homogeneous of degree $-1$, satisfy $\h(\X,\xi_i)=0$ and have the
property that $\{\X(u),\xi_1(u),\dots,\xi_n(u)\}$ is a basis for
$T_u\cq\subset T_u\aM$ for $u\in\cq$ close to $u_0$. Further, choose
an ambient vector field $Y$ such that $\h(\X,Y)=1$ (which forces $Y$
to be homogeneous of degree $-2$) and $\h(Y,\xi_i)|_{\cq}=0$ for all
$i$. Adding an appropriate multiple of $\X$, we may assume that
$\h(Y,Y)|_{\cq}=0$. Clearly, the fields $\X$, $\xi_i$ and $Y$ form a
frame for $T\aM$ in a neighbourhood of $u_0$ in $\cq$, and thus locally
around $u_0$. Let $\{\tilde Y,\eta_i,\tilde X\}$ be the local frame
dual to $\{\X,\xi_i,Y\}$. Then by construction $\tilde
Y|_{\cq}=Y|_{\cq}$ and $\tilde X|_{\cq}=\X|_{\cq}$, but this is not
true off $\cq$, since as we shall see immediately the weakest possible
assumption on Ricci flatness implies that $\h(\X,\X)$ is nonzero off
$\cq$.

\begin{thm*}
Let $\h$ be an ambient metric and let $\al$ be the one--form dual to
the infinitesimal generator $\X$ of the $\Bbb R_+$--action. Then we
have:

\noindent
(1) $\Ric(\h)(\X,\X)|_{\cq}=0$ if and only if
    $i_{\sX}d\al|_{\cq}=0$. If this is the case, then
    $r:=\tfrac{1}{2}\h(\X,\X)$ is a defining function for $\cq$ such
    that $Y\cdot r=1+O(r)$.

\noindent
(2) The tangential components of $\Ric(\h)(\X,\_)$ vanish along $\cq$
    if and only if $d\al|_{\cq}=0$. In particular, for any ambient
    metric $\h$ such that the tangential components of $\Ric(\h)$
    vanish along $\cq$, the procedure from \ref{2.4} can be used to
    obtain a normal standard tractor bundle.
\end{thm*}
\begin{proof}
From Proposition \ref{2.5} we get
$\h(\aR(\X,\zeta)\xi,\eta)=-\tfrac{1}{2}(\annd
d\al)(\zeta,\xi,\eta)$. Thus, we may compute $2\Ric(\h)(\X,\xi)$ by
taking the trace over $\zeta$ and $\eta$ in $(\annd
d\al)(\zeta,\xi,\eta)$, i.e.~by inserting the elements of dual frames
and summing up, and we use dual frames as introduced above. We only
have to consider the case that $\xi|_{\cq}$ is tangent to $\cq$, and
since we are only interested in the restriction of the result to
$\cq$, we may as well replace $\tilde X$ by $\X$ and $\tilde Y$ by
$Y$. The term with $\zeta=\X$ and $\eta=Y$ never contributes since
$(\annd d\al)(\X,\xi,Y)=-2\h(\aR(\X,\X)\xi,Y)=0$ by Proposition
\ref{2.5}.

Let us next look at the terms with $\zeta=\xi_i$ and $\eta=\eta_i$. By
definition tangential components of $\al$ vanish along $\cq$, so the
same holds for $d\al$. In particular, $d\al(\xi,\eta_i)$ vanishes
along $\cq$ and since $\xi_i$ is tangent to $\cq$ also $\xi_i\cdot
d\al(\xi,\eta_i)$ vanishes along $\cq$. Moreover, for any vector field
$\zeta$, the restriction of $\zeta-\al(\zeta)Y$ to $\cq$ is tangent to
$\cq$, which implies
$d\al(\xi,\zeta)|_{\cq}=\al(\zeta)d\al(\xi,Y)|_{\cq}$. Applying the
same argument with $\xi$ replaced by $\eta_i$, we see that
$d\al(\annd_{\xi_i}\xi,\eta_i)|_{\cq}=
\al(\annd_{\xi_i}\xi)d\al(Y,\eta_i)|_{\cq}$.  Since $\xi|_{\cq}$ is
tangent to $\cq$ and thus $\h(\xi,\X)|_{\cq}=0$, we see that, along
$\cq$, we have
$\al(\annd_{\xi_i}\xi)=\h(\annd_{\xi_i}\xi,\X)=-\h(\xi,\annd_{\xi_i}\X)$. 
Since both $\xi|_{\cq}$ and $\xi_i|_{\cq}$ are tangent to $\cq$ and
tangential components of $d\al$ vanish, formula \eqref{nX} from
\ref{2.3} implies that this restricts to $-\h(\xi,\xi_i)$ on $\cq$, so 
$$
d\al(\annd_{\xi_i}\xi,\eta_i)|_{\cq}=\h(\xi,\xi_i)d\al(\eta_i,Y). 
$$
Similarly,
$d\al(\xi,\annd_{\xi_i}\eta_i)=\al(\annd_{\xi_i}\eta_i)d\al(\xi,Y)$,
and $\al(\annd_{\xi_i}\eta_i)|_{\cq}=-\h(\eta_i,\xi_i)=-1$. Together
with the above, this implies that for any vector field $\xi$ such that
$\xi|_{\cq}$ is tangent to $\cq$, we obtain
$$
2\Ric(\h)(\X,\xi)|_{\cq}=(\annd d\al)(Y,\xi,\X)|_{\cq}
+nd\al(\xi,Y)|_{\cq}-\sum_i\h(\xi,\xi_i)d\al(\eta_i,Y)|_{\cq}  
$$

\noindent
(1) Putting $\xi=\X$ in the above formula, we see that the first
    summand vanishes since $\annd d\al$ is skew symmetric in the last
    two entries. On the other hand, the last sum vanishes since
    $\h(\X,\xi_i)=0$ by construction. Thus, we are left with
    $2\Ric(\h)(\X,\X)|_{\cq}=nd\al(\X,Y)|_{\cq}$, and since tangential
    components of $d\al$ vanish, the vanishing of $d\al(\X,Y)|_{\cq}$
    is equivalent to $i_{\X}d\al|_{\cq}=0$. We have already verified
    in \ref{2.3} that the latter condition implies that
    $r=\tfrac{1}{2}\h(\X,\X)$ is a defining function for $\cq$ since
    $dr|_{\cq}=\al|_{\cq}$. The last statement obviously implies
    $Y\cdot r=1+O(r)$. 

\noindent
(2) We may assume that the equivalent conditions of (1) are satisfied
    and show that vanishing of the rest of $\Ric(\h)(\X,\_)$ is
    equivalent to $d\al=O(r)$. Using the above formula for
    $2\Ric(\h)(\X,\xi)$, we first note that since
    $i_{\X}d\al$ vanishes along $\cq$, we get 
$$
(\annd d\al)(Y,\xi,\X)|_{\cq}=Y\cdot
    d\al(\xi,\X)|_{\cq}-d\al(\xi,\annd_Y\X)|_{\cq}. 
$$
    The first term in the right hand side may be written as $-Y\cdot
    (i_{\X}d\al(\xi))|_{\cq}$, and since $i_{\X}d\al|_{\cq}=0$ and
    $\xi|_{\cq}$ is tangent to $\cq$, this equals
    $di_{\X}d\al(\xi,Y)=\Cal L_{\X}d\al(\xi,Y)$. By construction,
    $\al$ is homogeneous of degree two, so the same holds for $d\al$,
    whence this gives $2d\al(\xi,Y)$. 
    For the second summand, we get
    $d\al(\xi,\annd_Y\X)|_{\cq}=\al(\annd_Y\X)d\al(\xi,Y)$ as above,
    and clearly $\h(\annd_Y\X,\X)=\tfrac{1}{2}Y\cdot\h(\X,\X)$, which
    restricts to $1$ on $\cq$ by part (1).

    Finally, by construction $\sum_i\h(\xi,\xi_i)\eta_i$ coincides
    with $\xi$ up to addition of a multiple of $\X$, so since
    $i_{\X}d\al|_{\cq}=0$ we obtain
    $\sum_i\h(\xi,\xi_i)d\al(\eta_i,Y)=d\al(\xi,Y)$. 

    Collecting our results, we see that (assuming
    $\Ric(\h)(\X,\X)|_{\cq}=0$) we get
    $2\Ric(\h)(\X,\xi)=nd\al(\xi,Y)$ for any $\xi$ such that
    $\xi|_{\cq}$ is tangent to $\cq$, which immediately implies the
    result.
\end{proof}

\section{An application}\label{3}
In this section, we show how our results can be applied to the study
of conformal invariants obtained from the ambient metric construction.
Some of these ideas were sketched in \cite{Gover-Peterson} but here we
generalise the setting considerably. In particular, we derive an
algorithm that can be used to compute a tractor expression for any
conformal invariant which can be obtained from the ambient metric
construction. Our results are however more general than that, since
they also deal with the case of ambient metrics which are Ricci flat
to higher order than those whose existence is proved by Fefferman and
Graham. While the existence of such metrics is obstructed on general
conformal manifolds, we believe studying these ``better'' metrics in
the cases when they do exist is very interesting.

It should be remarked at this point that another line of applications
of the results derived in this paper can be found in
\cite{Gover-Peterson}, where they are applied to the study of
conformally invariant powers of the Laplacian and $Q$--curvatures.

\subsection{}\label{3.1}
To carry out some computations, we introduce  abstract index
notation. Given an ambient manifold $\aM$ and an ambient metric $\h$
for a conformal structure $\cq\to M$, we write $\cce(w)$ for the space
smooth functions on $\aM$ which are homogeneous of degree $w$,
i.e. $\tf\in \cce(w)$ means $\X\cdot\tf =w\tf$. We will write $\cce^A
= \cce^A(0)$ ($\cce^A_{\cq}=\cce^A_{\cq}(0)$) to denote the space of
sections of $ T\aM$ ($ T\aM|_{\cq}$) which are homogeneous of degree
$-1$. (We adopt this convention since sections of $ \cce^A_{\cq}$
correspond to sections of the standard tractor bundle.) Then finally
we will write $\cce^{AB}(w)$ ($ \cce^{AB}_{\cq}(w)$) to mean
$\cce^A\otimes\cce^B\otimes \cce(w)$ ($ \cce^A_{\cq}\otimes
\cce^B_{\cq} \otimes \cce_{\cq}(w) $ respectively) and so forth. For
lower indices, the appropriate convention is that $\cce_A=\cce_A(0)$
denotes the space of ambient one--forms, which are homogeneous of
degree $1$, since then $\h_{AB}$ (which is homogeneous of degree $2$)
induces an isomorphism $\cce^A\to\cce_A$. The extensions to multiple
lower and mixed indices, as well as the notation for sections along
$\cq$ is done as above. In this context we will refer to $w$ as the
conformal {\em weight} (to distinguish it from the homogeneity
degree). This means that for an ambient tensor field, the conformal
weight equals the homogeneous degree plus the number of upper indices
minus the number of lower indices. We raise and lower indices using
the ambient metric $\h_{AB}$ and its inverse $\h^{AB}$. We also adopt
the usual conventions that round brackets (square brackets) around
indices indicate a symmetrisation (antisymmetrisation) of the enclosed
indices, except indices between vertical lines, and that the same
index occurring twice indicates a trace.  

We start with some general results about ambient metrics: 
\begin{prop*}
Let $\pi:\cq\to M$ be a conformal structure, $\aM$ an ambient manifold
and $\h$ an ambient metric on $\aM$ with curvature
$\aR=\aR_{AB}{}^C{}_D$ and Ricci curvature
$\aRic=\aRic_{AB}=\aR_{CA}{}^C{}_B$. Then we have:

\noindent
(1) $\annd^E\aR_{EABC}=2\annd_{[B}\aRic_{C]A}$.

\noindent
(2) $\aDe\aR_{ABCD}=2(\annd_A\annd_{[C}\aRic_{D]B}-
\annd_B\annd_{[C}\aRic_{D]A})+\Ps_{ABCD}$, where $\aDe$ denotes the
ambient Laplacian and $\Ps_{ABCD}$ is a linear combination of partial
contractions of $\aR\otimes\aR$. 

\noindent
(3) Let $\Ph_{A\dots B}\in\ce_{A\dots B}(w)$ be any section. Then the
    commutator of the Laplacian $\aDe$ with a covariant derivative
    $\annd$ acts as
\begin{align*}
[\aDe,\annd_C]\Ph_{A\dots B}=&-2\h^{EF}\left(\Ph_{E\dots
B}\annd_{[F}\aRic_{A]C}+\dots+\Ph_{A\dots E}\annd_{[F}\aRic_{B]C}\right)\\
&-2\left(\aR_{EC}{}^F{}_A\annd^E\Ph_{F\dots B}+\dots+
\aR_{EC}{}^F{}_B\annd^E\Ph_{A\dots F}\right)\\
&+\aRic_{CE}\annd^E\Ph_{A\dots B},
\end{align*}
where in the two sums there is one summand for each index of $\Ph$,
and $I$ is contracted into that index. 
\end{prop*}
\begin{proof}
(1) The algebraic Bianchi identity $0=\aR_{[EA|B|C]}$ together with
the usual symmetries of the Riemann curvature gives us
$\aR_{EABC}=-\aR_{ABEC}+\aR_{ACEB}$. The differential Bianchi identity
$0=\annd_{[F}\aR_{AB]EC}$ together with the symmetries of $\aR$ leads
to $\annd_F\aR_{ABEC}=-\annd_B\aR_{FAEC}+\annd_A\aR_{FBEC}$, and
similarly we get
$\annd_F\aR_{ACEB}=-\annd_C\aR_{FAEB}+\annd_A\aR_{FCEB}$. Contracting
with $\h^{EF}$ the claim now follows from symmetry of $\aRic$. 

\noindent
(2) By definition
    $(\aDe\aR)_{ABCD}=\annd^E\annd_E\aR_{ABCD}$. Using the
    differential Bianchi identity and curvature symmetries, we may
    write $\annd_E\aR_{ABCD}$ as
    $-\annd_B\aR_{EACD}+\annd_A\aR_{EBCD}$. Now the commutator of two
    covariant derivatives is given by the algebraic action of the
    curvature, so $-\annd^E\annd_B\aR_{EACD}$ may be written as the sum
    of $-\annd_B\annd^E\aR_{EACD}$ and a sum of partial contractions of
    $\aR\otimes\aR$. Similarly,
    $-\annd^E\annd_A\aR_{EBCD}$ is the sum of
    $-\annd_A\annd^E\aR_{EBCD}$ and a sum of partial contractions of
    $\aR\otimes\aR$. Now the result immediately follows from (1).

\noindent
(3) Let us compute $\aDe\annd_C\Ph_{A\dots
    B}=\h^{EF}\annd_E\annd_F\annd_C\Ph_{A\dots B}$. The definition of
    the curvature reads as $[\annd_A,\annd_B]V^C=\aR_{AB}{}^C{}_DV^D$,
    and thus $[\annd_A,\annd_B]V_C=-\aR_{AB}{}^D{}_CV_D$. Using this,
    we get 
$$
\annd_F\annd_C\Ph_{A\dots B}=\annd_C\annd_F\Ph_{A\dots
B}-(\aR_{FC}{}^I{}_A\Ph_{I\dots B}+\dots+ \aR_{FC}{}^I{}_B\Ph_{A\dots
I}),
$$
with one summand for each index of $\Ph$ in the sum in
brackets. Hitting that sum with $\annd^F$, each summand splits into a
sum of two terms, one in which $\annd^F$ acts on $\aR$ and one in
which $\annd^F$ acts on $\Ph$. Using (1) we see that the terms in
which $\annd^F$ acts on $\aR$ exactly give the terms in the first sum
of the claimed formula for $[\aDe,\annd_C]\Ph_{A\dots B}$. On the
other hand, the terms in which $\annd^F$ acts on $\Ph$ exactly give
half of the second sum in the claimed formula.

Again swapping covariant derivatives, we may write
$\annd_E\annd_C\annd_F\Ph_{A\dots B}$ as the sum of
$\annd_C\annd_E\annd_F\Ph_{A\dots B}$ (which after contraction with
$\h^{EF}$ gives $\annd_C\aDe\Ph_{A\dots B}$) and
$$
-\aR_{EC}{}^I{}_F\annd_I\Ph_{A\dots B}-(\aR_{EC}{}^I{}_A\annd_F\Ph_{I\dots
B}+\dots+\aR_{EC}{}^I{}_B\annd_F\Ph_{A\dots I}),
$$
again with one summand for each index of $\Ph$ in the sum in
brackets. Contracting with $\h^{EF}$, the the sum in brackets gives
the second half of the second sum in our claimed formula, while the
other summand gives the last term in the claimed formula.
\end{proof}

\subsection*{Remark} 
Of course, in the proof of part (2), it is no problem to
compute an explicit formula for the sum of partial contractions
$\Ps_{ABCD}$ of $\aR\otimes\aR$ (see \cite{Gover-Peterson}). 

\subsection{}\label{3.2}
To proceed, we next specialise to an ambient metric $\h$ such that
tangential components of the ambient Ricci curvature $\aRic$ vanish
along $\cq$. By part (2) of Theorem \ref{2.6} the one form $\al$ dual
to the infinitesimal generator $\X$ of the $\Bbb R_+$--action then
satisfies $d\al=O(r)$ and the procedure of \ref{2.4} can be applied to
construct a normal standard tractor bundle $(\Cal T,h,\nabla)$ from
$(\aM,\h)$.

We may regard the ambient curvature $\aR$ as 2-form taking values in
$\End(T\aM)$. We have observed in \ref{2.5} above that if $\tilde\xi$
and $\tilde\eta$ are invariant lifts of vector fields
$\xi,\eta\in\frak X(M)$, then $\aR(\tilde\xi,\tilde\eta)$ is precisely
the homogeneous degree 0 section of $\End(T\aM|_{\cq})$ corresponding
to the section $R(\xi,\eta)$ of $\End{\Cal T}$.
Thus the homogeneous $\End(T\aM)$
valued three form $\alpha\wedge \aR$ is, along $\cq$, uniquely
determined by the tractor curvature. Similarly, the Levi-Civita
connection $\annd$ is determined by its action on vector fields
homogeneous of degree $-1$, so along $\cq$ covariant derivatives in
tangential directions are determined by the underlying tractor
connection.

Since $d\al|_{\cq}=0$, we know from \ref{2.3} that
$r:=\tfrac{1}{2}\h(\X,\X)$ satisfies $\al=dr+O(r)$, and hence $r$ is a
smooth defining function for $\cq$. Notice that the ambient vector
field $\X=\X^A$ has conformal weight $1$ and since the ambient
covariant derivative is compatible with homogeneities, the ambient
differential operator $\annd_A$ has conformal weight $-1$.  By
definition, the ambient one--form $\alpha_A$ is given by
$\h_{AB}\X^B$, so we will also denote this form by $\X_A$. (So $\X$
will mean either a 1-form or a vector field according to index
placement and/or context.) For example $\alpha\wedge \aR$ is
$3\X_{[A}\aR_{BC]}{}^D{}_E\in\cce_{[ABC]}{}^D{}_E (-1)$. Note that by
definition $\X^A\X_A=\h_{AB}\X^A\X^B=2r$.

To compute efficiently in the sequel, we have to determine the
commutators of covariant derivatives with $r$ and $X_A$, viewed as
multiplication operators. Since $d\al=O(r)$, we get
$i_{\X}d\al=2(\al-dr)=O(r)$, so $\al=dr+r\be$ for some ambient
one--form $\be$. Moreover, $d\al=dr\wedge\be+rd\be$ and this being
$O(r)$ implies vanishing of the tangential components of $\be$ along
$\cq$, whence $\be=\ph \al+r\tilde\ga$ for some ambient smooth
function $\ph$ and one--form $\tilde\ga$, and hence $\al=dr+\ph
r\al+r^2\tilde\ga$. In particular, $dr=(1-\ph r)\al+O(r^2)$, and
viewing $r$ as a multiplication operator, this implies the commutation
formula $[\annd_A,r]=(1-\ph r)\X_A+O(r^2)$. Further, the above
equations immediately imply $dr\wedge\al=O(r^2)$ and
$d\al=r\ga\wedge\al+O(r^2)$, where $\ga:=d\ph-2\tilde\ga$, which in
index notaiton reads as $\annd_{[A}\X_{B]}=r\ga_{[A}\X_{B]}$. On the
other hand, equation \eqref{nX} from \ref{2.3} gives us
$\annd_A\X_B=\h_{AB}+r\ga_{[A}\X_{B]}+O(r^2)$.  Viewing $X_B$ as a
multiplication operator, we thus get the commutator formula
$[\annd_A,\X_B]=\h_{AB}+r\ga_{[A}\X_{B]}+O(r^2)$.

The next step is to compute two basic tractor operators in ambient
terms. The first obvious operator to consider is $\alpha\wedge \annd$,
which, along $ \cq$, obviously only needs derivatives in tangential
directions, and may be written as
$\adD_{AB}:=2\X_{[A}\annd_{B]}$. Since $[\annd_A,r]=(1-\ph
r)\X_A+O(r^2)$, one immediately concludes that $[\adD_{AB},r]=O(r^2)$,
which in particular means that for any ambient tensor field $V$, the
restriction $(\adD_{AB} V)|_{\cq}$ depends only on $V|_{\cq}$. Hence
for arbitrary indices, we get a well defined operator
$\adD_{AB}:(\cce_{\cq})^{C\dots D}_{E\dots
F}(w)\to(\cce_{\cq})^{C\dots D}_{[AB]E\dots F}(w)$, that clearly can
be computed in terms of the underlying standard tractor bundle. It is
easy to identify this operator: By definition, the adjoint tractor
bundle of $M$ is the bundle of orthogonal endomorphisms of the
standard tractor bundle $ \Cal T$. Recall sections of $ \Cal T$ may be
identified with sections in $\cce_{\cq}^A$ and vice versa. Using the
inverse $\h^{AB}$ of the ambient metric, there is similarly a
one-to-one correspondence between smooth sections of the adjoint
tractor bundle $\Cal A$ and sections in $\cce_{\cq}^{[AB]}$.  Thus
$\adD_{AB}$ determines a conformally invariant operator $\dD $ on $ M$
which goes between ${\Cal A}\otimes {\Cal F}\otimes \ce[w] $ and
${\Cal F}\otimes \ce[w] $, where ${\Cal F}$ is any tensor power of
$\Cal T$.  There is a natural projection $\Cal A\to TM$. Under the
identification of $\Ga(\Cal A)$ with $\cce_{\cq}^{[AB]}$ this is
explicitly given by mapping $\Ph^{AB}$ to $\X_A\Ph^{AB}-\X_A\Ph^{BA}$
modulo multiples of $\X^B$.  Using this, one immediately verifies that
on the standard tractor bundle, $\dD$ coincides with the composition
of the tractor connection with this projection.  On density bundles
one obtains a similar composition of Levi-Civita connection with the
projection invariantly combined with the canonical action of an
adjoint tractor on the density bundle; in total the fundamental
$D$--operator, see \cite[section 3]{luminy}. The obvious compatibility
of $\adD_{AB}$ with tensor products then implies that it is exactly
the operator obtained by twisting the fundamental $D$ on the density
bundle with the tractor connection on the tractor bundle. This is
precisely the ``double-D'' operator of \cite{gosrni,conf-invar} (and
see also \cite[section 3]{luminy}).
 
Now we can follow these sources to obtain the tractor
$D$--operator: Consider the operator $\h^{AB}\adD_{A(Q}\adD_{|B|P)_0}$,
which by construction acts tangentially on tensor fields along
$\cq$. Here $(\cdots)_0$ indicates the trace-free
symmetrisation over enclosed indices (excluding any in the $|\cdots|$).
Using the commutator formulae from above, one immediately verifies
directly that 
\begin{align*}
4\h^{AB}\X_{[A}\annd_{Q]}\X_{[B}\annd_{P]}=&-n\X_Q\annd_P-
\h_{PQ}\X^A\annd_A+\X_Q\X_P\aDe\\
&-\X_Q\X^A\annd_A\annd_P-\X_P\X^A\annd_Q\annd_A+O(r).
\end{align*}
Next, one easily verifies that on any homogeneous tensor field, the
operator $\X^A\annd_A\annd_P+\X^A\annd_P\annd_A$ acts in the same way
as $-2\annd_P+2\annd_P\X^A\annd_A$. Thus we conclude that
$\h^{AB}\adD_{A(Q}\adD_{|B|P)_0}=-\X_{(Q}\D_{P)_0}+O(r)$ where
$\D_A:=(n-2)\annd_A +2\annd_A \X^B\annd_B -\X_A \aDe$.  Since the map
$\xi_P\mapsto \X_{(Q}\xi_{P)_0}$ is an injection of $T^*\aM|_\cq$ into
$\otimes^2T^*\aM|_\cq$ we see that this construction determines $\D_A
$ as a well defined operator between $ \cce_\cq^A\otimes
\cce_\cq^\psi(w)$ and $ \cce_\cq^\Psi(w-1)$, where $ \cce^\Psi$ is any
tensor power of $ \cce^B$. Thus $ \D$ determines an operator $ D$
between weighted tractor bundles on $ M$.  By its construction from $
\adD$ it is clear that this action of $\D_A$ is determined by the
underlying standard tractor bundle and its connection.  In fact this
construction of $ \D$ is exactly the interpretation on $\cq $ of the
construction of the standard tractor D operator from $\dD$ as in
\cite{gosrni} and \cite[3.2]{luminy}. For which it follows immediately
that $ D$ is the standard tractor $D$ operator. 

We can easily verify explicitly that $\D_A$ acts tangentially on
homogeneous ambient tensor fields along $\cq$.  On ambient tensors of
conformal weight $w$ we have
$\D_A=(n+2w-2)\annd_A-\X_A\aDe$. Moreover, from above we know that
$\annd_A r=(1-\ph r)X_A+O(r^2)$. Hitting this with $\annd^A$, we
immediately conclude that $\aDe r=(n+2)+O(r)$. More generally, if $\V$
is any ambient tensor field of conformal weight $w$, then
$\aDe(rV)=(\aDe r)V+2\h^{AB}(\annd_A r)\annd_B V+r\aDe
V=(n+2w+2)V+O(r)$. Since $r\V$ has conformal weight $w+2$,
together with the above formula for the action of $\D_A$ on tensor
fields of fixed conformal weight this implies that $\D_A rV=O(r)$, so
$\D_A$ indeed acts tangentially along $\cq$.

\subsection{}\label{3.3}
The simple relation between the operator $\D_A$ and the ambient
covariant derivative together with the fact that $\D_A$ acts
tangentially leads to very remarkable consequences. The point here is
that since, along $ \cq$, $\D_A$ depends only on the underlying
standard tractor bundle, and for any ambient tensor field $\Ph$ the
restriction of $\D_A\Ph$ to $\cq$ depends only on the restriction of
$\Ph$ to $\cq$. In particular, we may apply this to the ambient
curvature $\aR$ and its covariant derivatives. The restriction of
$\aR$ to $\cq$ may be viewed as a section of
$(\cce_{ABCD})_{\cq}(-2)$, and similarly for any $\ell>0$ the
restriction of $\annd^\ell\aR=(\annd\o\dots\o\annd)\aR$ to $\cq$
defines a section of $(\cce_{A_1\dots A_{\ell+4}})_{\cq}(-2-\ell)$.
\begin{prop*}
Let $\aM$ be an ambient manifold for a conformal structure $\cq\to M$
on an $n$--dimensional manifold $M$, and let $\h$ be an ambient metric
on $\aM$ with curvature $\aR$ and Ricci curvature $\aRic$. If for some
$k>0$ we have $\aRic_{AB}=O(r^{k+1})$ and tangential components of
$\aRic$ vanish to order $k+2$ along $\cq$, and if $n\neq 2k+4$, then
there is a universal formula that computes $\annd^k\aR|_{\cq}$ from
$\aR|_{\cq}$ and $\annd^\ell\aR|_{\cq}$ for $\ell<k$.
\end{prop*}
\begin{proof}
The section $\annd^{k-1}\aR$ has conformal weight $-k-1$, so from above
we know that
$\D_A\annd^{k-1}\aR=(n-2k-4)\annd_A\annd^{k-1}\aR-\X_A\aDe\annd^{k-1}\aR$,
and by assumption $n-2k-4\neq 0$. Since $\D_A\annd^{k-1}\aR|_{\cq}$
depends only on $\annd^{k-1}\aR|_{\cq}$ it suffices to compute
$\aDe\annd^{k-1}\aR$ from the restrictions of the sections
$\annd^\ell\aR$ for $\ell<k$ to $\cq$. Now one immediately verifies
inductively that
$$
[\aDe,\annd^{k-1}]=\sum_{\ell=0}^{k-2}\annd^\ell[\aDe,\annd]\annd^{k-2-\ell}.
$$
Our assumptions on $\aRic$ together with part (3) of proposition
\ref{3.1} imply that for any ambient tensor field $\Ph_{A\dots B}$, we
can write $[\aDe,\annd_C]\Ph_{A\dots B}$ as 
$$
-2\left(\aR_{EC}{}^F{}_A\annd^E\Ph_{F\dots B}+\dots+
\aR_{EC}{}^F{}_B\annd^E\Ph_{A\dots F}\right)+O(r^k).
$$
Inserting $\Ph=\annd^{k-2-\ell}\aR$ for some $\ell=0,\dots,k-2$, we
get an expression for $[\aDe,\annd]\annd^{k-2-\ell}\aR$ in terms of
$\aR$ and $\annd^{k-1-\ell}\aR$ up to some $O(r^k)$. Applying
$\annd^\ell$ and restricting to $\cq$, the $O(r^k)$ cannot contribute,
and we only get covariant derivatives of order at most $k-1$ of
$\aR$. Thus we obtain a universal formula which expresses 
$[\aDe,\annd^{k-1}]\aR|_{\cq}$ in terms of $\aR|_{\cq}$ and
$\annd^\ell\aR|_{\cq}$ for $\ell<k$. 

To complete the proof it hence suffices to analyse
$\annd^{k-1}\aDe\aR$, and we have computed $\aDe\aR$ in part (2) of
Proposition \ref{3.1}. Now we claim that the term
$\annd_A\annd_{[C}\aRic_{D]B}- \annd_B\annd_{[C}\aRic_{D]A}$ showing
up in that formula is $O(r^k)$. Indeed, by assumption
$\aRic=O(r^{k+1})$ and tangential components of $\aRic$ vanish to
order $k+2$ along $\cq$, which implies that
$\aRic_{AB}=r^{k+1}(\X_AK_B+\X_BK_A)+O(r^{k+2})$ for some ambient
one--form $K_B$. Using the commutation formula for $\annd$ and $r$, we
see that $\annd_C\aRic_{DB}=(k+1)r^k\X_C(\X_DK_B+\X_BK_D)+O(r^{k+1})$,
so skewing over $C$ and $D$, we are left with
$(k+1)r^k\X_BK_{[D}\X_{C]}+O(r^{k+1})$. Hitting this with $\annd_A$ we
get $k(k+1)r^{k-1}\X_A\X_BK_{[D}\X_{C]}+O(r^k)$, and skewing over $A$
and $B$ the claim follows. 

But then it follows from part (2) of Proposition \ref{3.1} that
$\annd^{k-1}\aDe\aR|_{\cq}=\annd^{k-1}\Ps_{ABCD}|_{\cq}$, and
since $\Ps_{ABCD}$ is a partial contraction of $\aR\otimes\aR$ we
conclude that $\annd^{k-1}\Ps_{ABCD}|_{\cq}$ can be expressed by a
universal formula in terms of $\aR|_{\cq}$ and
$\annd^\ell\aR|_{\cq}$ for $\ell\leq k-1$. 
\end{proof}

\subsection{}\label{3.4}
We have noted in \ref{3.2} that the restriction of the section
$\X_{[A}\aR_{BC]EF}$ of $\cce_{[ABC][EF]}(-1)$ to $\cq$ depends only
on the tangential components of the ambient curvature, which equal
with the curvature of the normal standard tractor
connection. Consequently, $3\D^A\X_{[A}\aR_{BC]EF}$ is a section of
$\cce_{[BC][EF]}(-2)$ whose restriction to $\cq$ can be computed from
the tractor curvature (and thus depends only on the underlying
conformal structure). To compute this explicitly, we need the
commutator of the Laplacian with $X_A$. From section 3.2 above, we
have $\annd_B X_A=\h_{BA}+r\ga_{[B}\X_{A]}+X_A\annd_B+O(r^2)$. Hitting
this with $\annd^B$, we obtain $\aDe X_A=2\annd_A+\tfrac{1}{2}
X^B\ga_B X_A+X_A\aDe+O(r)$.  Since $X^B\ga_B $ and the curvature
$\aR_{BCEF}$ both have conformal weight $-2$, this allows us to write
$\D_G\X_{[A}\aR_{BC]EF}$ as $(n-2)\annd_G\X_{[A}\aR_{BC]EF}+
(\frac{1}{2}\X^J\ga_J-\aDe)\X_G\X_{[A}\aR_{BC]EF}+O(r)$.  Taking into
account the conformal weights and using the formulae derived above, we
can now expand this expression explicitly and contracting with
$\h^{GA}$ we obtain
\begin{equation}\label{curvid}
\begin{array}{lll}
\lefteqn{3\D^A\X_{[A}\aR_{BC]EF}=}&&
\\
&&
(n-2)[(n-4)\aR_{BCEF}+2\X_{[B}\annd^A\aR_{C]AEF}]
\\
&&
+( X^J\ga_J -2\aDe)(\X^A\X_{[B}\aR_{C]AEF})+O(r). 
\end{array}
\end{equation}

In particular, in dimensions $\neq 4$, the curvature of the ambient
metric shows up in this formula. In these dimensions, existence of an
ambient metric $\h$ such that $\Ric(\h)=O(r)$ and tangential
components of $\Ric(\h)$ vanish to second order along $\cq$ has been
proved in \cite{F-G}. So let us assume that we deal with such a
metric.  We first claim that the Ricci condition implies
$d\al=O(r^2)$. We know that $d\al=O(r)$, so $d\al=r\be$ for some
two--form $\be$, which is homogeneous of degree $0$, and thus has
conformal weight $-2$. From \ref{3.2} we conclude that
$\aDe(d\al)=\aDe(r\be)=(n-2)\be+O(r)$, so it suffices to prove
$\aDe(d\al)=O(r)$ in order to conclude $d\al=O(r^2)$. In index
notation, proposition \ref{2.5} reads as
$\annd_A(d\al)_{BC}=2\X^E\aR_{EABC}$.  Consequently, we may compute
$\aDe(d\al)$ as $\tfrac{1}{2}\annd_A\X_E\aR^{EA}{}_{BC}$. From
\ref{3.2} we know that up to an $O(r)$ we may replace $\annd_A\X_E$ by
$\h_{AE}+\X_E\annd_A$, so the fact that the curvature is skew in the
first two indices implies
$\aDe(d\al)=\tfrac{1}{2}\X^E\annd^A\aR_{EABC}+O(r)$. By part (1) of
proposition \ref{3.1}, this may be rewritten as
$-\X^E\annd_{[B}\aRic_{C]E}$. As in the proof of Proposition
\ref{3.3}, we may write $\aRic_{CE}=r(\X_CK_E+\X_EK_C)+O(r^2)$ for an
appropriate $K$ and applying the argument from that proof with $k=0$,
we see that $\annd_{[B}\aRic_{C]E}=\X_EK_{[C}\X_{B]}+O(r)$.
Contracting this with $\X^E$, we get an $O(r)$ term, which completes
the proof that $d\al=O(r^2)$.

From above we know that
$\X^A\aR_{CAEF}=-\tfrac{1}{2}\annd_C(d\al)_{EF}$, and tensoring with
$\X_B$ and skewing over $B$ and $C$, we see that 
$\X^A\X_{[B}\aR_{C]AEF}$ equals $-\tfrac{1}{4}\adD_{BC}(d\al)_{EF}$. Since
$d\al=O(r^2)$, the same is true for $\adD(d\al)$, and from \ref{3.2}
we know that applying $\aDe$ to an $O(r^2)$ term we get an $O(r)$
term, so the whole last line in \eqref{curvid} restricts to zero on
$\cq$. On the other hand, using part (1) of Proposition \ref{3.1} we
get
$$
\X_{[B}\annd^A\aR_{C]AEF}=(-\X_B\annd_{[E}\aRic_{F]C}+
\X_C\annd_{[E}\aRic_{F]B})
$$
As above, we may write $\aRic_{FC}=r(X_FK_C+X_CK_F)+O(r^2)$ and
then $\annd_{[E}\aRic_{F]C}=X_CK_{[F}X_{E]}+O(r)$. Multiplying by
$X_B$ and antisymmetrising in $B$ and $C$ we obtain
$\X_{[B}\annd^A\aR_{C]AEF}=O(r)$. Collecting our results, we see that
assuming that assuming that $\Ric(\h)=O(r)$ and tangential components
of $\Ric(\h)$ vanish to second order along $\cq$, formula
\eqref{curvid} boils down to
$3\D^A\X_{[A}\aR_{BC]EF}|_{\cq}=(n-2)(n-4)\aR_{BCEF}|_{\cq}$. Using
this, we now get

\begin{thm*}
Let $\aM$ be an ambient manifold for a conformal structure $\cq\to M$
on an $n$--dimensional manifold $M$, and let $\h$ be an ambient metric
on $\aM$ with curvature $\aR$ and Ricci curvature $\aRic$. Assume that
$\Ric(\h)=O(r^{k+1})$ and tangential components of $\Ric(\h)$ vanish
to order $k+2$ along $\cq$ for some $k\geq 0$. Then using the
convention that $\annd^0\aR=\aR$ we have:

\noindent
(1) If $n$ is odd or $n$ is even and $k<\frac{n-4}{2}$, then for each
    $0\leq\ell\leq k$ there is a universal tractor formula that
    computes $\annd^\ell\aR$ from the tractor curvature of the
    underlying standard tractor bundle.

\noindent
(2) If $n$ is even and $\frac{n-4}{2}<\ell\leq k$, then there is a
    universal tractor formula that computes $\annd^\ell\aR$ from the
    tractor curvature of the underlying standard tractor bundle and
    $(\annd^{\frac{n-4}{2}}\aR)|_{\cq}$.
\end{thm*}
\begin{proof}
Since our assumptions on $\Ric(\h)$ imply that at least
$\Ric(\h)=O(r)$ and tangential components of $\Ric(\h)$ vanish to
second order, we may can compute $\aR$ via
$\aR_{BCEF}|_{\cq}=\frac{3}{(n-2)(n-4)}\D^A\X_{[A}\aR_{BC]EF}|_{\cq}$
provided that $n\neq 4$. Thus we get (1) for $\ell=0$. Iterated
application of Proposition \ref{3.3} then leads to a universal formula
for $\annd^\ell\aR$ in terms of $\aR$ provided that in no step we get
$n=2\ell+4$, and so (1) follows. 

For (2), we first note that by (1) we get a universal formula for
$\annd^i\aR$ in terms of the tractor curvature for
$i<\frac{n-4}{2}$. But using this, the result again follows from
iterated application of Proposition \ref{3.3}. 
\end{proof}

\subsection*{Remarks}
Part (1) of this theorem ties in nicely with the results on existence
and uniqueness of ambient metrics in \cite{F-G}. It is shown in that
paper that in odd dimensions there exists an infinite order power
series solution along $\cq$ for an ambient metric $\h$ such that
$\Ric(\h)$ vanishes to infinite order along $\cq$ and this solution is
unique up to the action of an equivariant diffeomorphism fixing $\cq$.
On the other hand, in the case of even dimensions the Fefferman-Graham
construction is obstructed at finite order. More precisely, if $n$ is
even then there exist an ambient metric $\h$ such that
$\Ric(\h)=O(r^{\frac{n-4}{2}})$ and tangential components of
$\Ric(\h)$ vanish to order $\frac{n-4}{2}+1$, but there is an
obstruction against the existence of an ambient metric such that
$\Ric(\h)=O(r^{\frac{n-4}{2}+1})$. Moreover, this solution is unique
up to the action of an equivariant diffeomorphism fixing $\cq$ and
addition of terms which vanish to order $\frac{n-2}{2}+1$ along $\cq$.
In particular, this implies that for $n\neq 4$ the curvature of such
an ambient metric as well as those of its covariant derivatives that
are covered in part (1) of the theorem are intrinsic to the underlying
conformal structure. Since the tractor curvature is intrinsic to the
underlying conformal structure, part (1) of the theorem provides an
alternative proof for this fact, which clearly comes close to an
alternative proof for uniqueness of Fefferman--Graham metrics. It
seems to us, however, that the ideas developed in this paper should
also have applications to the question of existence of ambient metrics
and the nature of the obstruction in the Fefferman--Graham
construction, so we will take up this whole circle of problems
elsewhere. Below we will show how our results can be applied to obtain
explicit descriptions of ambient Weyl invariants.

On the other hand, part (2) of the theorem goes significantly beyond
the results in \cite{F-G}, since it analyses the cases in which the
obstruction in the Fefferman--Graham construction vanishes. It shows
that in these cases the only essential new ingredient is the critical
covariant derivative $\annd^{\frac{n-4}{2}}\aR$ of the ambient
curvature, which then determines all higher covariant derivatives by
universal tractor formulae. It should also be remarked here that large
parts of this critical covariant derivative are again determined by
the underlying conformal structure, since covariant derivatives in
tangential directions are determined by the tractor connection.

\subsection{Applications to the study of ambient Weyl
  invariants}\label{3.5} 
One of the main applications of the conformal ambient metric
construction is that it allows a systematic construction of conformal
invariants. It is well known that by Weyl's classical invariant theory
any Riemannian invariant can be written as a linear combination of
complete contractions of tensor powers of iterated covariant
derivatives of the Riemann curvature tensor. Consider an arbitrary
Riemannian invariant in odd dimensions or an invariant depending only
on covariant derivatives up to order less than $\tfrac{n-4}{2}$ in
even dimensions $n\neq 4$. Applied to a Fefferman--Graham metric, one
can consider the restriction of the resulting function to $\cq$. Since
parts of different homogeneity of this function must be individually
invariant, we may without loss of generality assume that our function
is homogeneous of some degree and hence may be interpreted as a
density on $M$. Now the fact that we deal with a Riemannian invariant
exactly eliminates the diffeomorphism freedom in the ambient metric,
while the freedom of adding terms that vanish along $\cq$ to the
appropriate order has already been taken care of. Consequently, the
resulting density on $M$ is conformally invariant. Invariants obtained
in that way are called ambient Weyl--invariants. 

Theorem \ref{3.4} not only provides an alternative proof for the fact
the the construction outlined above leads to conformal invariants, but
also provides a way to compute explicit formulae for ambient Weyl
invariants, which otherwise is a difficult problem. On the one hand,
Theorem \ref{3.4} directly leads to an iterative way to compute
tractor expressions for ambient Weyl invariants. For many purposes
this is already sufficient, since one obtains genuine formulae for the
given invariant and many qualitative features of the invariant can be
appreciated in this compact form. On the other hand, converting
tractor expressions into formulae in terms of a metric representing
the conformal class, its Levi-Civita connection and curvature is a
completely mechanical procedure which even may be left to a computer. 
Since we do not want to introduce too much tractor calculus here, we
only roughly analyse a simple case below. More involved applications
in a similar direction can be found in \cite{Gover-Peterson}. 

For the first step in this analysis, we assume $n\neq 4$, and we are
dealing with an ambient metric $\h$ such that $\Ric(\h)=O(r)$ and
tangential components of $\Ric(\h)$ vanish to second order. Then we
have observed above that
$3\D^A\X_{[A}\aR_{BC]EF}|_{\cq}=(n-2)(n-4)\aR_{BCEF}|_{\cq}$. Further,
$\X_{[A}\aR_{BC]EF}|_{\cq}$ depends only on the tangential components
(in the indices $B$ and $C$) of the ambient curvature which are given
by the tractor connection. Otherwise put, extending the tractor
curvature $\ka_{bcEF}$ in any way to a tractor field $\Om_{BCEF}$ and
forming $\X_{[A}\Om_{BC]EF}$, one will obtain the same result. In the
general setting of tractor calculus, the expression $W_{BC}{}^E{}_F:=
\frac{3}{n-2}D^AX_{[A} \Omega_{BC]}{}^E{}_F$ has been known for some
time (see \cite{gosrni,conf-invar}) as a natural conformally invariant
tractor extension of the tractor curvature (which itself is an
extension of the Weyl curvature or the Cotton--York tensor in
dimension $3$). Thus we see that $W_{BC}{}^E{}_F$ is the tractor
field equivalent to $(n-4)\aR_{BC}{}^E{}_F|_\cq$.

To describe a formula for $W_{BC}{}^E{}_F$ we have to introduce some
basic elements of tractor calculus, see \cite{BEG, gosrni, conf-invar,
Gover-Peterson}. We use lower case letter for tensor indices and upper
case letters for standard tractor indices and also ambient indices. We
use the tractor metric and its inverse to raise and lower tractor
indices. Choosing a metric $g$ from the conformal class, we may raise
and lower tensor indices using $g$ (but taking into account that this
changes the weight), and the standard tractor bundle $\ce^A$ splits as
$\ce[1]\oplus\ce_a[1]\oplus\ce[-1]$. This can be conveniently encoded
by adding to the natural section $X^A\in\ce^A[1]$ (which represents
the natural inclusion $\ce[-1]\to\ce^A$) tractor sections
$Y^A\in\ce^A[-1]$ and $Z^{aA}\in\ce^{aA}[-1]$ which represent the
other two inclusions that depend on the choice of the metric
$g$. Basic properties of the tractor metric imply that $Y^AX_A=1$,
$X^AX_A=Y^AY_A=X^AZ^a_A=Y^AZ^a_A=0$ and
$Z^{aA}Z^b{}_A=g^{ab}$. Further, we denote by $W_{abcd}$ the
Weyl--curvature, by $S$ the scalar curvature of $g$, by
$\Rho_{ab}=\frac{1}{n-2}(\Ric_{ab}-\frac{1}{2(n-1)}Sg_{ab})$ the
rho--tensor of $g$ and by $B_{eb}:=\nd^q\nd^pW_{peqb}
+(n-3)\Rho^{qp}W_{peqb}$ the Bach tensor. Using the formulae in
\cite{gosrni, conf-invar} one verifies that $W_{ABCE}$ is given by
\begin{align*}
&(n-4) Z_A{}^aZ_B{}^bZ_C{}^cZ_E{}^e
W_{abce} -4(n-4)Z_A{}^aZ_B{}^bX_{[C}Z_{E]}{}^e \nd_{[a}\Rho_{b]e} \\ 
-4&(n-4)X_{[A}Z_{B]}{}^b Z_C{}^cZ_E{}^e \nd_{[c}\Rho_{e]b}
+\frac{4}{n-3} X_{[A}Z_{B]}{}^b X_{[C} Z_{E]}{}^e B_{eb}
\end{align*}

This is a general tractor formula not related to the ambient metric in
any way, so in particular, it also holds in dimension $n=4$. In that
case, only the last term survives, which shows that $B_{ab}$ is a
conformal invariant in dimension $4$. However, our interest here is in
the case $n\neq 4$, and the main fact that we need from the above
formula is that since $X^AZ^a_A=0$ and $Z^{aA}Z^b{}_A=g^{ab}$ one
immediately concludes that any complete contraction of
$W_{ABCD}\otimes\dots\otimes W_{IJKL}$ gives the complete contraction
of $W_{abcd}\otimes\dots\otimes W_{ijkl}$ with the same pairing of
indices. Applying this to the case of ambient Weyl invariants, we
obtain an alternative proof of \cite[Proposition 3.2]{F-G} that any
complete contraction of a tensor power of the ambient curvature gives
the ``same'' complete contraction of the corresponding tensor power of
the Weyl--curvature. Of course, these are exactly the invariants one
knew about in advance, so to get more interesting invariants one has
to consider covariant derivatives of the ambient curvature.

The simplest ambient Weyl invariant involving covariant derivatives is
$\|\annd \aR\|^2:=(\annd^A \aR^{BCEF})\annd_A \aR_{BCEF}$. In order
that this is well defined, we have to assume that $n\neq 4,6$ and that
we are dealing with an ambient metric $\h$ such that $\Ric(\h)=O(r^2)$
and tangential components of $\Ric(\h)$ vanish to third order along
$\cq$. To apply our method, we first note following the proof of
Proposition \ref{3.3} that $\aR$ has conformal weight $-2$, and thus 
$$
\D_A\aR_{BCEF}=(n-6)\annd_A\aR_{BCEF}-\X_A\aDe\aR_{BCEF}. 
$$
From part (2) of Proposition \ref{3.1} we next conclude that our
assumptions on $\Ric(\h)$ imply that the restriction to $\cq$ of
$\aDe\aR_{BCEF}$ coincides with a partial contraction of
$\aR\otimes\aR$. Thus we see that up to a nonzero factor and the
addition of complete contractions of tensor powers of the Weyl
curvature, the ambient Weyl invariant $\|\annd\aR\|^2$ coincides with 
$(D^AW^{BCEF})D_AW_{BCEF}$. Using the standard formulae for the
tractor D--operator and the formulae for $W_{BCEF}$ above, this can be
expanded into formula in terms a chosen representative of the
conformal class. While this is straightforward, it is quite tedious
since there are many components in $D_AW_{BCEF}$ and most of them
do not contribute to the final invariant.

A more efficient way to proceed is to rewrite the original invariant
as follows: Applying the differential Bianchi identity to the first
term, we see that $(\annd_A \aR_{BCEF})\annd^A \aR^{BCEF}$ may be
written as $2(\annd_B \aR_{ACEF})\annd^A \aR^{BCEF}$. On the other
hand, consider $\annd^A\annd_B\aR_{ACEF}$. Switching the covariant
derivatives can be compensated by adding partial contractions of
$\aR\otimes\aR$, and using part (1) of Proposition \ref{3.1} and the
proof of Proposition \ref{3.3}, we see that $\annd^A\aR_{ACEF}$ can be
written as $4r\X_CK_{[E}\X_{F]}+O(r^2)$ for an appropriate tractor
field $K_E$. Hitting this with $\annd_B$, we obtain
$4\X_B\X_CK_{[E}\X_{F]}+O(r)$, which vanishes upon contraction into
$\aR^{BCEF}$. The upshot of this is that, along $\cq$ and up to adding
complete contractions of tensor powers of $\aR$, we may rewrite
$\|\annd\aR\|^2$ as $2\annd_A\annd_B \aR^{ACDE}\aR^B{}_{CDE}$. 

Following the same idea as for the single covariant derivative above, one
next shows that
$$
\D_A\D_B\aR^{ACDE}\aR^B{}_{CDE}=(n-6)(n-8)\annd_A\annd_B
\aR^{ACDE}\aR^B{}_{CDE},
$$
up to complete contractions of tensor powers of $\aR$. Thus the
conformal invariant given by $||\annd\aR||^2 $ is, up to the addition
of complete contractions of tensor powers of the Weyl curvature,
exactly
$$
\frac{2}{(n-4)^2(n-6)(n-8)}D_AD_BW^{ACDE}W^B{}_{CDE}.
$$
in dimensions other than $4,6$. Note that we know in advance that
the expression $D_AD_BW^{ACDE}W^B{}_{CDE}$ will yield $ (n-8)$ times a
conformal invariant because $||\annd\aR||^2 $ is well defined in
dimension 8.  This is then easily expanded out using formulae for the
tractor connection as in \cite{luminy} or \cite{Gover-Peterson} to
yield (once again modulo complete contractions of tensor powers of
$W$),
$$
\frac{2}{n-6}\square \| W\|^2 +\frac{n-10}{(n-6)}[\|\nd W\|^2 + 8h(W,U)
- 4(n-6)\|C\|^2],
$$
where $\square$ is the tractor Laplacian, $C_{abd}=2\nd_{[a}\Rho_{b]d}$
is the Cotton--York tensor and
$U_{abcd}:=\nd_aC_{cdb}+\Rho_a{}^eW_{ebcd}$. This shows explicitly how
the invariant simplifies in dimension 10. Alternatively we can
re-express in the more compact form
$$
\|\nd W\|^2+16 (W,U) -4(n-10) \|C\|^2.
$$
It is readily verified that this precisely agrees with the result
obtained by Fefferman-Graham in \cite{F-G}. (In fact we have borrowed
some notation from that source to simplify the comparison.) Note that
although $\|\annd \aR\|^2 $ is not well defined in dimensions 4 and 6,
the last display does give an invariant in these dimensions.  This is
immediate from the fact that in arbitrary dimension $n$,
$D_AD_BW^{ACDE}W^B{}_{CDE}$ may be written as the sum of
$\frac {(n-4)^2(n-6)(n-8)}{2} (||\nd W||^2+16 h(W,U) -4(n-10) ||C||^2)$
and complete contractions of tensor powers of $W$, and that the
transformation formulae for the scalars $\|\nd W\|^2$, $h(W,U)$ and
$\|C\|^2$, under a change of metric from the conformal class, is of 
polynomial type.


\begin{thebibliography}{XX}

\bibitem{BEG} T.N. Bailey, M.G. Eastwood, A.R. Gover, \textit{Thomas's
structure bundle for conformal, projective and related structures},
Rocky Mountain J. \textbf{24} (1994), 1191--1217.

\bibitem{nonlocal} T.\ Branson, A.R.\ Gover, \textit{Conformally
    Invariant Non-Local Operators}, Pacific J.\ Math.\ {\bf 201}
  (2001) 19--60.

\bibitem{Cartan} E. Cartan, \textit{Les espaces \`a connexion
conforme}, Ann. Soc. Pol. Math., \textbf{2} (1923), 171--202.

\bibitem{luminy} A. \v Cap, A.R. Gover, \textit{Tractor bundles for
irreducible parabolic geometries}, SMF S\'eminaires et congr\`es
{\bf 4} (2000) 129--154, electronically available at
http://smf.emath.fr/SansMenu/Publications/SeminairesCongres/

\bibitem{tractors} A. \v Cap, A.R. Gover, \textit{Tractor Calculi for
Parabolic Geometries}, Trans. Amer. Math. Soc. \textbf{354} (2002),
1511-1548, electronically available as Preprint ESI 792 at
http://www.esi.ac.at

\bibitem{Cap-Schichl} A. \v Cap, H. Schichl, \textit{Parabolic
Geometries and Canonical Cartan Connections}, Hokkaido
Math. J. \textbf{29} No.3 (2000) 453-505

\bibitem{Mike-Srni} M.G. Eastwood, \textit{Notes on Conformal
Differential Geometry}, Supp. Rend.  Circ.  Matem.  Palermo,
\textbf{43} (1996), 57--76.

\bibitem{F-G} C. Fefferman and C.R. Graham, \textit{Conformal
invariants}, in \'{E}lie Cartan et les Math\'{e}matiques
d'Adjourd'hui, (Ast\'{e}risque, hors serie), (1985), 95--116.

\bibitem{GrF} C.\ Fefferman and C.R.\ Graham,  {\em $ Q$-curvature
    and Poincar\'{e} metrics},  preprint. math.DG/0110271.

\bibitem{gosrni} A.R. Gover, \textit{Aspects of parabolic invariant theory}, 
 Supp.\ Rend.\ Circ.\ Matem.\ Palermo, Ser.\ II, Suppl.\ {\bf 59}
(1999) 25--47.

\bibitem{conf-invar} A.R. Gover, \textit{Invariants and calculus for
    conformal geometry}, Adv.\ Math. {\bf 163} (2001), 206--257.

\bibitem{Gover-Peterson}A.R. Gover, L.J. Peterson, \textit{Conformally
invariant powers of the Laplacian, Q-curvature, and tractor calculus},
42 pp., Preprint math-ph/0201030, electronically available at
http://arXiv.org. 

\bibitem{Graham} C.R. Graham, \textit{Conformally invariant powers of
the Laplacian, II: Nonexistence} J. London Math. Soc., \textbf{46}
(1992), 566--576.

\bibitem{GJMS} C.R. Graham, R. Jenne, L.J. Mason, G.A. Sparling, 
\textit{Conformally invariant powers of the Laplacian, I: Existence}, J. London
Math. Soc., \textbf{46} (1992), 557--565

\bibitem{GrZ} C.R. Graham and M.\ Zworski, {\em Scattering matrix in conformal 
geometry},  preprint. math.DG/0109089 

\bibitem{GrW} C.R. Graham and E. Witten, {\em Conformal anomaly of
    submanifold observables in AdS/CFT correspondence}, Nuclear Phys.\ 
  B, \textbf{546} (1999) 52--64, hep-th/9901021.

\bibitem{HS2} M. Henningson and K. Skenderis, {\em Holography and the
Weyl anomaly, Proceedings of the 32nd International Symposium
Ahrenshoop on the Theory of Elementary Particles (Buckow, 1998)},
Fortschr.\ Phys.\ \textbf{48} (2000) 125--128, hep-th/9812032.

\bibitem{T} T.Y. Thomas, {\em On conformal geometry},
  Proc.\ Natl.\ Acad.\ Sci.\ USA
  {\bf 12} (1926) 352--359.

\bibitem{Thomasbook} T.\ Y.\ Thomas, ``The Differential Invariants of
Generalized Spaces,'' Cambridge University Press , Cambridge, 1934. 

\bibitem{W} E.\ Witten,  {\em Anti-de Sitter space and
    holography}, Adv.\ Theor.\ Math. Phys.\, {\bf 2} (1998) 253--291,
hep-th/9802150.

\end{thebibliography}
\end{document}